\documentclass[12pt, english]{article}
\usepackage[utf8]{inputenc}
\usepackage{amsthm}
\usepackage{amssymb}
\usepackage[all]{xy}
\usepackage{amsmath}
\usepackage{tikz-cd}
\usepackage{tikz}
\usepackage{float}
\usepackage{bbm}

\tikzcdset{scale cd/.style={every label/.append style={scale=#1},
		cells={nodes={scale=#1}}}}
\newcommand{\RNum}[1]{\uppercase\expandafter{\romannumeral #1\relax}}

\usepackage{rotating}
\usepackage{caption, booktabs}
\usepackage{verbatim}
\usetikzlibrary{calc}

\usepackage{tikz}
\usetikzlibrary{patterns}

\usepackage{pdflscape}

\usepackage{mathtools}
\usepackage[english]{babel}
\usepackage{enumitem}

\usepackage{hyperref}
\usepackage{bookmark}

\newtheorem{definition}{Definition}[section]

\newtheorem{proposition}[definition]{Proposition}
\newtheorem{remark}[definition]{Remark}
\newtheorem{lemma}[definition]{Lemma}

\newtheorem{theorem}[definition]{Theorem}
\newtheorem{corollar}[definition]{Corollary}
\newtheorem{construction}[definition]{Construction}

\newtheorem{example}[definition]{Example}

\DeclareMathOperator{\ord}{ord}
\DeclareMathOperator{\Hom}{Hom}
\DeclareMathOperator{\Pic}{Pic}
\DeclareMathOperator{\coker}{coker}
\DeclareMathOperator{\Lie}{Lie}
\DeclareMathOperator{\Aut}{Aut}

\DeclareMathOperator{\Gr}{Gr}

\DeclareMathOperator{\Stab}{Stab}
\DeclareMathOperator{\Int}{Int}
\DeclareMathOperator{\Min}{Min}
\DeclareMathOperator{\Supp}{Supp}
\DeclareMathOperator{\Bound}{Bound}
\DeclareMathOperator{\Cone}{Cone}
\DeclareMathOperator{\Convhull}{Convhull}

\begin{document}

	\title{Jacobian rings and the infinitesimal Torelli Theorem}
\author{Julius Giesler,\\ University of T\"ubingen}
\date{\today}
\maketitle
\begin{abstract}
	In this article we deal with jacobian rings and identify a mixed Hodge component of a nondegenerate hypersurface in the torus with a lattice geometric quotient vector space. We introduce a period map, study its differential and compute the kernel of the differential much explicitly via certain Laurent polynomials. As a main application we deal with the infinitesimal Torelli theorem (ITT). We study the kernel of the cohomological map for explicit deformations and complete our study by dealing with the remaining part $\coker(\kappa_{\mathbb{P},f})$ (cokernel of the Kodaira-Spencer map) in dimensions $n \geq 4$.
\end{abstract}

\section{Introduction}

Being a subsequent paper to the paper (\cite{Gie22}) summarizing, generalizing and extending results from my doctoral thesis (\cite{Gie23}) we adopt the notation respectively. Given a \textit{nondegenerate} Laurent polynomial 
\begin{align} \label{represenation_f_intro}
	f = \sum\limits_{m \in \Delta \cap M} a_m \cdot x^m
\end{align}
with $n$-dimensional Newton polytope $\Delta$ take 
\[ Z_f := \{ f=0 \} \subset T. \]

We are interested in the mixed Hodge components 
	\begin{align} \label{equation_mixed_Hodge_components_of_Z_f_of_min_wei}
		H^{n-k,k-1} H^{n-1}(Z_f, \mathbb{C}) \cong R_{Int,f}^{k}
	\end{align}
	of weight $n-1$ and their dimensions, where $R_{Int,f}^k$ denotes a graded component of an \textit{interior module} defined via jacobian rings (of Batyrev), (see \cite{Bat93}). The isomorphism (\ref{equation_mixed_Hodge_components_of_Z_f_of_min_wei}) identifies the mixed Hodge components of minimal weight with vector spaces defined by the lattice geometry of integral multiples of $\Delta$ (definition \ref{definition_R_Int_f_k}). In this article we find a representation of $R_{Int,f}^k$ as quotient vector space that illustrates Griffiths concept of \textit{(infinitesimal) variation of Hodge structure} for toric hypersurfaces: (Proposition \ref{proposition_Hodge_comp_Z_f_arbitary_k_dimensions_and_basis})

\begin{proposition} \label{proposition_Hodge_comp_Z_f_arbitary_k_dimensions_and_basis_intro}
	Let $\Delta$ be an $n$-dimensional lattice polytope with $l^*(\Delta) > 0$ and a given $f$. Given $\Gamma_1,...,\Gamma_{n+1} \leq \Delta$ with $n_{\Gamma_1},...,n_{\Gamma_{n+1}}$ affine linear independent. Then 
	\begin{align*}
		L^*(k \cdot \Delta) / U_{f,k} = R_{Int,f}^k \quad k=1,...,n
	\end{align*}
	where $U_{f,k}$ denotes the vector space over $\mathbb{C}$ spanned by
	\begin{align}
		& g_{\Gamma_i}(f) \cdot x^v \quad i=1,...,n+1, \quad v \in \Int((k-1) \cdot \Delta) \cap M \label{formula_relation_g_Gamma_i_f_the_1}	 \\
		& g_{\Gamma}(f) \cdot x^v \quad \Gamma \leq \Delta, \quad \quad \quad \quad \, \, \, \, v \in \Int((k-1) \cdot \Gamma) \cap M \label{formula_relation_g_Gamma_i_f_the_2}
	\end{align}
	with certain Laurent polynomials $g_{\Gamma}(f)$. For $k=2$ these relations are linear independent. Here $\Int$ denotes the relative interior and $L^*(k \cdot \Delta)$ the $\mathbb{C}$-vector space with basis $x^m$, where $m \in \Int(k \cdot \Delta) \cap M$.
\end{proposition}

The $g_{\Gamma}(f)$'s in fact are useful generators of the first component $J_{\Delta,f}^1$ of the jacobian ideal. 
Using a resolution of the vector space $R_{Int,f}^k$ we were able to specify a complete list of relations between the above generators for $k > 2$.
\leavevmode
\\


We then introduce the usual period map $\mathcal{P}_{B,f}$ with entries
\begin{align*}
	\mathcal{P}_{B,f}^k: f' \mapsto F^{k}H^{n-1}(Y_{f'},\mathbb{C}) \subset H^{n-1}(Y_f,\mathbb{C})
\end{align*}
with $\{ F^k \}$ the Hodge filtration. This map is known to be holomorphic. The $k$-th component $d \mathcal{P}_{B,f}^k$ of the differential $d \mathcal{P}_{B,f}$ factors as follows
\begin{equation}
	\begin{tikzcd} \label{Commutative_diagram_Kodaira_Spencer_map_dP_intro}
		T_{B,f}  \arrow[swap]{dr}{d \mathcal{P}_{B,f}^k } \arrow{r}{\kappa_{f}} & H^{1}(Y_f,T_{Y_f})  \arrow{d}{\Phi_f^k} \\
		& \Hom(H^{k}(Y_f, \Omega_{Y_f}^{n-1-k}), H^{k+1}(Y_f, \Omega_{Y_f}^{n-2-k}))
	\end{tikzcd}
\end{equation}
Our idea is strainghtforward: Compute the kernel $\ker(d \mathcal{P}_{B,f})$ and compare it with $\ker(\kappa_f)$, which we explicitly computed in the preceeding article \cite{Gie22}, to finally deduce on $\ker(\Phi_{f \vert{Im(\kappa_f)}})$. A detailed study of the components $R_{Int,f}^k$ and the relations between the generators of $0 \in R_{Int,f}^k$ gives the following key result (theorem \ref{Theorem_formula_ker_ITT})

\begin{theorem} \label{Theorem_formula_ker_ITT_intro}
	Let $\Delta$ be an $n$-dimensional lattice polytope with $(0,...,0) \in \Int(\Delta) \cap M$. 
	If $f \in L(\Delta)$ is nondegenerate
	\begin{align} \label{general_formula_for_ker_DP_MDelta_handable_intro}
		\ker(d \mathcal{P}_{B,f}^1) \overset{\mod \Lie(T)}{=}&  \Big\langle g_{\Gamma}(f) \cdot x^w \in L(\Delta)/\mathbb{C} \cdot f \quad \mid \quad  \Gamma \leq \Delta \textrm{ a facet}, \nonumber  \\
		& w + v \in \big( \Int(\Delta) \cup \Int(\Gamma) \big) \cap M, \nonumber \quad \\ 
		& \forall \, v \in \Int(\Delta) \cap M \Big\rangle. \nonumber
	\end{align}
	where $\langle \rangle$ denotes the span as $\mathbb{C}$-vector space.
\end{theorem}

This direct approach to the diagram of Griffiths (\ref{Commutative_diagram_Kodaira_Spencer_map_dP_intro}) is not published yet and yields a both computational and much explicit perspective. 
\leavevmode 
\\
Let $g = g_{\Gamma}(f) \cdot x^w \in \ker(d \mathcal{P}_{B,f}^1)$. Then we get the distinction
\begin{equation} \label{distinction_scalar_product_0_minus1_getminus2}
	\begin{tikzcd}
		\langle w,n_{\Gamma} \rangle = \left\{\begin{array}{lll} 0, & \Leftrightarrow g \equiv 0 \\
			-1, &  \Leftrightarrow g \in \ker(\kappa_f) \\
			\leq -2, & \quad \, \, \textrm{exceptional cases}
		\end{array}\right.
	\end{tikzcd}
\end{equation}
and $\langle w,n_{\Gamma'} \rangle \geq 0$ for $\Gamma' \neq \Gamma$ by Theorem (\ref{Theorem_formula_ker_ITT_intro}) with $n_{\Gamma}$ the inner facet normal to $\Gamma$. The last case in the distinction causes $\Phi_{f \vert{Im \, \kappa_f}}$ to be not injective, which implies that the infinitesimal Torelli theorem (ITT) fails for such hypersurfaces $Y_f$. 
\begin{corollar}
	Given $\Int(\Delta) \cap M \nsubseteq \{\textrm{plane}\}$ and $f \in L(\Delta)$ nondegenerate then $\ker(d \mathcal{P}_{B,f}) = \ker(\kappa_f)$. Thus
\[  \ker(\Phi_{f \vert{Im \, \kappa_f}}) = \{ 0  \},  \]
We say that $Y_f$ fulfills the ITT for $\mathcal{X} \rightarrow B$.  
\end{corollar}

This corollary follows easily from above: Namely we construct a $3$-dimensional empty polytope $Q$ with $6$ vertices. $Q$ is the convex span of an empty triangle with vertex $(0,0,0)$ and this triangle dilated by $w$. Then we use a theorem of White on empty $3$-simplices to deduce $|Q \cap M| \geq 7$ under the assumption $\langle w,n_{\Gamma} \rangle \leq -2$, a contradiction.
\leavevmode 
\\ 

Finally we tackle the homomorphism $\Phi_{f \vert{\coker \, \kappa_f}}$ in dimension $n \geq 4$: Given a \textit{functorial toric resolution} of singularities $\tilde{\mathbb{P}} \rightarrow \mathbb{P}$ we have
\begin{align*}
	H^0(\tilde{\mathbb{P}}, T_{\tilde{\mathbb{P}}}) \cong H^0(\mathbb{P}, T_{\mathbb{P}}).
\end{align*}
We write $\tilde{Y}_f$, or by abuse of notation also $Y$, for the smooth closure of $Z_f$ in $\tilde{\mathbb{P}}$. The results on the respective \textit{Kodaira Spencer map} $\kappa_{\tilde{\mathbb{P}},f}$ remain valid (section \ref{section_Kod_Sp_map_smooth_bir_model}): Up to slightly adapting $\Delta$ we get $\kappa_{\tilde{\mathbb{P}},f} = \kappa_f$ and given $n \geq 4$ we have an exact sequence
\begin{align*}
	0 \rightarrow Im(\kappa_{f}) \rightarrow H^1(\tilde{Y}, T_{\tilde{Y}}) \rightarrow H^1(\tilde{\mathbb{P}}, T_{\tilde{\mathbb{P}}}) \rightarrow 0.
\end{align*}
The cohomological map $\Phi_f$ in the infinitesimal Torelli theorem (ITT) is made up of cup product and contraction. Given $n \geq 4$ we show how $\Phi_f$ could be realized on $\coker(\kappa_{f})$ again by cup product and contraction (section \ref{section_ITT_on_coker_computation}). We then deduce on the ITT on the remaining part $\coker(\kappa_{f})$ (section \ref{section_ITT_on_coker_computation_2}) and arrive at:

\begin{theorem}
	Let $\Delta \subset M_{\mathbb{R}}$ be an $n$-dimensional lattice polytope with
	\begin{align*}
		n \geq 4, \quad \quad \Int(\Delta) \cap M \nsubseteq \textrm{ plane}
	\end{align*}
	and $f \in U_{reg}(\Delta)$ (with Newton polytope $\Delta$). Then the smooth birational model $Y = \tilde{Y}_f$ of $Z_f$ from above fulfills the infinitesimal Torelli theorem (ITT), that is
	\begin{align*}
		\Phi_f = \bigoplus_p \Phi_f^p: H^1(Y,T_Y) \rightarrow \bigoplus_{p} \Hom \Big( H^{n-1-p}(Y, \Omega_Y^p), H^{n-p}(Y, \Omega_Y^{p-1}) \Big)
	\end{align*}
	is injective. 
\end{theorem}

Our article relies on results from (\cite{DK86}) and (\cite{Bat93}). Let us list common methods with the authors V.I. Danilov, A.G. Khovanskii and V.V. Batyrev:
\begin{itemize}
	\item We reduce the calculation of the Hodge component $H^{n-k,k-1}(Y_f, \mathbb{C})$ to the calculation of (\ref{equation_mixed_Hodge_components_of_Z_f_of_min_wei}) (similarly to \cite{DK86}), see section \ref{Section_cohomology_affine_hypersurface_nonsingular_compactification}.
	\item We follow \cite{Bat93} to define the \textit{jacobian ring} $R_f$ and the \textit{(interior) module} $R_{Int,f}$ over $R_f$ (see section \ref{section_constr_jac_ring}) and of course use the identification (\ref{equation_mixed_Hodge_components_of_Z_f_of_min_wei}).
\end{itemize}

In section (\ref{section_smooth_and_stable_points}) we study $\Delta$ up to replacing by an integral lattice vector. In particular we effort $(0,...,0) \in \Int(\Delta)$  to build a smooth quotient of $B$ by the torus $T$.

\section{The Jacobian ring of Batyrev} \label{section_constr_jac_ring}


\begin{definition}
	Let $\Delta$ be an $n$-dimensional lattice polytope and $S_{\Delta}$ denote the subalgebra of $\mathbb{C}[x_{0},x_{1}^{\pm},...,x_{n}^{\pm}]$  spanned as $\mathbb{C}$-vector space by $1$ and all monomials 
	\[  x_{0}^{k}x_1^{m_1}...x_{n}^{m_{n}}, \quad \textrm{where } k \in \mathbb{N}_{\geq 1},  \quad m_1,...,m_n \in \mathbb{Z}, \]
	such that the rational point
	\[ \frac{m}{k} := \Big( \frac{m_{1}}{k},...,\frac{m_{n}}{k} \Big) \]
	belongs to $\Delta$. Denote the $k$-th graded piece of $S_{\Delta}$ by $S_{\Delta}^k$.	The subalgebra $S_{\Delta}^*$ is defined in the same way except that we require that $m/k$ belongs to the interior of $\Delta$.
\end{definition}

We identify lattice points $m \in M$ with monomials $x^m \in \mathbb{C}[x_1^{\pm 1},...,x_n^{\pm 1}]$.

\begin{construction} \label{construction_Laurent_polynomials_F_i} \cite{Bat93} \\
	\normalfont
	Given a Laurent polynomial $f$ in the variables $x_1,...,x_n$ with $n$-dimensional Newton polytope $\Delta$ let
	\[ F(x_{0},x_{1},...,x_{n}) := x_{0}f(x_1,...,x_n) - 1  \]
	which is an equation for the complement $T \setminus Z_{f} \subset \overline{T}:= (\mathbb{C}^{*})^{n+1}$ of $Z_f=\{f=0\} \subset T$. Consider the logarithmic derivatives 
	\[ F_{i}(x_{0},x) := x_{i}\frac{\partial}{\partial x_{i}} F(x_{0},x) = x_i \frac{\partial (x_0 f)}{\partial x_i} \quad 0 \leq i \leq n. \]
\end{construction}
\begin{definition} \cite{Bat93} \\
	\normalfont
	The graded ideal $J_{\Delta,f}$ within $S_{\Delta}$ generated by $F_{0},...,F_{n}$ is called the \textit{jacobian ideal} and the quotient ring
	\[ R_{f} := S_{\Delta}/J_{\Delta,f} \]
	is called the \textit{Jacobian ring of Batyrev} (associated to the Laurent polynomial $f$). We denote the $k$-th homogeneous component of $R_{f}$ by $R_{f}^k$.
\end{definition}

\begin{theorem} \label{theorem_characterization_Delta_regular} \cite[Thm.4.8]{Bat93} \\
	The jacobian ring $R_{f}$ is a graded ring. It is finite dimensional as $\mathbb{C}$-vector space if and only if $f$ is nondegenerate. In this case the dimensions of $R_f^k$ are independent of the polynomial $f$. $f$ is nondegenerate if and only if $F_0,...,F_n$ are algebraically independent over $\mathbb{C}$.
\end{theorem}

\begin{definition} \label{definition_R_Int_f_k} (\cite[Prop.9.0.6]{Bat93}) \\
	We denote the homogeneous ideal $R_f \cap S_{\Delta}^*$ of $R_{f}$ by $R_{Int,f}$ and its k-th homogeneous component by $R_{Int,f}^k$. We call $R_{Int,f}$ the \textit{interior $R_{f}$-module}. Let $D_f^k := L^*(k \cdot \Delta)/(F_0,...,F_n) \cdot L^*((k-1) \cdot \Delta)$.
\end{definition}

\begin{remark} \label{remark_duality_Poincare_duality}
	\normalfont
	Poincaré duality on $H^{n-1}(Y,\mathbb{C})$ restricts to an isomorphism
	\[ (R_{Int,f}^k)^* \cong R_{Int,f}^{n+1-k}, \quad \textrm{and } (R_{\Delta,f}^k)^* \cong D_f^{n+1-k}  \]
	by (\cite[Remark 9.5, Prop.9.7]{Bat93}). 
\end{remark}


\section{Construction of the components $R_{Int,f}^k$} \label{section_examples_comp_jac_ring_interor_module}

\begin{remark}
	\normalfont	
	By definition
	\[ R_{Int,f}^1 \cong L^*(\Delta) \] 
	(independently of $f$). More interestingly
	\begin{align} \label{formula_RDelta_Int_2_relations_Jf_Delta}
		R_{Int,f}^k \cong L^{*}(k \Delta) \big/ \, \Big( J_{\Delta,f}^k \cap L^{*}(k \Delta) \Big).
	\end{align}
	Since $J_{\Delta,f}^k$ is a graded ideal of $S_{\Delta}$ inductively
	\[ J_{\Delta,f}^k = L((k-1)\Delta) \cdot J_{\Delta,f}^1  \]
	\underline{Aim:} Switching to different generators of $J_{\Delta,f}^1$ allows us to describe the relations in $R_{Int,f}^k$ more explicitly.
	
\end{remark}

\begin{figure}[H]
	
	\begin{center}
		
		\begin{tikzpicture}[scale=0.7]

			\begin{scope}[xshift = 10cm]
				
				\fill (-4,0) circle(2pt);
				\fill (-4,1) circle(2pt);
				\fill (-4,2) circle(2pt);
				\fill (-4,3) circle(2pt);
				\fill (-4,4) circle(2pt);
				\fill (-4,5) circle(2pt);
				\fill (-4,6) circle(2pt);
				
				\fill (-3,0) circle(2pt);
				\fill (-3,1) circle(2pt);
				\fill (-3,2) circle(2pt);
				\fill (-3,3) circle(2pt);
				\fill (-3,4) circle(2pt);
				\fill (-3,5) circle(2pt);
				\fill (-3,6) circle(2pt);
				
				\fill (-2,0) circle(2pt);
				\fill (-2,1) circle(2pt);
				\fill (-2,2) circle(2pt);
				\fill (-2,3) circle(2pt);
				\fill (-2,4) circle(2pt);
				\fill (-2,5) circle(2pt);
				\fill (-2,6) circle(2pt);
				
				\fill (-1,0) circle(2pt);	
				\fill (-1,1) circle(2pt);
				\fill (-1,2) circle(2pt);
				\fill (-1,3) circle(2pt);
				\fill (-1,4) circle(2pt);
				\fill (-1,5) circle(2pt);
				\fill (-1,6) circle(2pt);
				
				\fill (0,0) circle(2pt);
				\fill (0,1) circle(2pt);
				\fill (0,2) circle(2pt);
				\fill (0,3) circle(2pt);
				\fill (0,4) circle(2pt);
				\fill (0,5) circle(2pt);
				\fill (0,6) circle(2pt);
				
				\fill (1,0) circle(2pt);
				\fill (1,1) circle(2pt);
				\fill (1,2) circle(2pt);
				\fill (1,3) circle(2pt);
				\fill (1,4) circle(2pt);
				\fill (1,5) circle(2pt);
				\fill (1,6) circle(2pt);
				
				\fill (2,0) circle(2pt);
				\fill (2,1) circle(2pt);
				\fill (2,2) circle(2pt);
				\fill (2,3) circle(2pt);
				\fill (2,4) circle(2pt);
				\fill (2,5) circle(2pt);
				\fill (2,6) circle(2pt);
				
				\fill (3,0) circle(2pt);
				\fill (3,1) circle(2pt);
				\fill (3,2) circle(2pt);
				\fill (3,3) circle(2pt);
				\fill (3,4) circle(2pt);
				\fill (3,5) circle(2pt);
				\fill (3,6) circle(2pt);
				
				\fill (4,0) circle(2pt);
				\fill (4,1) circle(2pt);
				\fill (4,2) circle(2pt);
				\fill (4,3) circle(2pt);
				\fill (4,4) circle(2pt);
				\fill (4,5) circle(2pt);
				\fill (4,6) circle(2pt);
				
				\fill (5,0) circle(2pt);
				\fill (5,1) circle(2pt);
				\fill (5,2) circle(2pt);
				\fill (5,3) circle(2pt);
				\fill (5,4) circle(2pt);
				\fill (5,5) circle(2pt);
				\fill (5,6) circle(2pt);
				
				\fill (6,0) circle(2pt);
				\fill (6,1) circle(2pt);
				\fill (6,2) circle(2pt);
				\fill (6,3) circle(2pt);
				\fill (6,4) circle(2pt);
				\fill (6,5) circle(2pt);
				\fill (6,6) circle(2pt);
				
				\draw (-4,0) -- (2,6);
				\draw (2,6) -- (6,0);
				\draw (6,0) -- (-4,0);
				
				\node[left] at (-1.1,4.8) {{\bf $R_{Int,f}^2$}};
				\node[left] at (-3.0,2.5) {{\bf $v_1 + \Delta$}};
				\node[right] at (6.35,2.6) {{\bf $v_2+\Delta$}};
				\node[right] at (4.35,5.6) {{\bf $v_3+\Delta$}};
				
				\node[left] at (-0.9,1) {{\bf $v_1$}};
				\node[right] at (3.9,1) {{\bf $v_2$}};
				\node[above] at (2,3.9) {{ \bf $v_3$}};
				\node[left] at (-3.9,0) {{\bf $2v_1$}};
				\node[right] at (5.9,0) {{\bf $2v_2$}};
				\node[above] at (2,5.9) {{\bf $2v_3$}};
				
				\draw[->] (-2,4.5)--(1,2.5);
				\draw[->] (-3.2,2.3) -- (-1.7,1.7);
				\draw[->] (6.5,2.6) -- (4.3,1.6);
				\draw[->] (4.5, 5.6) -- (2.4,4.7);
				
				\draw (-1,1) -- (2,4);
				\draw (2,4) -- (4,1);
				\draw (4,1) -- (-1,1);
				
				\draw[dashed] (-1,3) -- (1,0);
				\draw[dashed] (-1.2,2.8) -- (2/3,0);
				\draw[dashed] (-1.4,2.6) -- (1/3,0);
				\draw[dashed] (-1.6,2.4) -- (0/3,0);
				\draw[dashed] (-1.8,2.2) -- (-1/3,0);
				\draw[dashed] (-2.0,2.0) -- (-2/3,0);
				\draw[dashed] (-2.2,1.8) -- (-3/3,0);
				\draw[dashed] (-2.4,1.6) -- (-4/3,0);
				\draw[dashed] (-2.6,1.4) -- (-5/3,0);
				\draw[dashed] (-2.8,1.2) -- (-6/3,0);
				\draw[dashed] (-3.0,1.0) -- (-7/3,0);
				\draw[dashed] (-3.2,0.8) -- (-8/3,0);
				\draw[dashed] (-3.4,0.6) -- (-9/3,0);
				\draw[dashed] (-3.6,0.4) -- (-10/3,0);
				\draw[dashed] (-3.8,0.2) -- (-11/3,0);
				
				\draw[dashed] (-1,3) -- (4,3);
				\draw[dashed] (-0.8,3.2) -- (4-2/15,3.2);
				\draw[dashed] (-0.6,3.4) -- (4-4/15,3.4);
				\draw[dashed] (-0.4,3.6) -- (4-6/15,3.6);
				\draw[dashed] (-0.2,3.8) -- (4-8/15,3.8);
				\draw[dashed] (-0.0,4.0) -- (4-10/15,4.0);
				\draw[dashed] (0.2,4.2) -- (4-12/15,4.2);
				\draw[dashed] (0.4,4.4) -- (4-14/15,4.4);
				\draw[dashed] (0.6,4.6) -- (4-16/15,4.6);
				\draw[dashed] (0.8,4.8) -- (4-18/15,4.8);
				\draw[dashed] (1.0,5.0) -- (4-20/15,5.0);	
				\draw[dashed] (1.2,5.2) -- (4-22/15,5.2);
				\draw[dashed] (1.4,5.4) -- (4-24/15,5.4);
				\draw[dashed] (1.6,5.6) -- (4-26/15,5.6);
				\draw[dashed] (1.8,5.8) -- (4-28/15,5.8);
				
				\draw[dashed] (1.0,0) -- (4,3);
				\draw[dashed] (1+5/15,0) -- (4+2/15,2.8);
				\draw[dashed] (1+2*5/15,0) -- (4+4/15,2.6);
				\draw[dashed] (1+3*5/15,0) -- (4+6/15,2.4);
				\draw[dashed] (1+4*5/15,0) -- (4+8/15,2.2);
				\draw[dashed] (1+5*5/15,0) -- (4+10/15,2.0);
				\draw[dashed] (1+6*5/15,0) -- (4+12/15,1.8);
				\draw[dashed] (1+7*5/15,0) -- (4+14/15,1.6);
				\draw[dashed] (1+8*5/15,0) -- (4+16/15,1.4);
				\draw[dashed] (1+9*5/15,0) -- (4+18/15,1.2);
				\draw[dashed] (1+10*5/15,0) -- (4+20/15,1.0);
				\draw[dashed] (1+11*5/15,0) -- (4+22/15,0.8);
				\draw[dashed] (1+12*5/15,0) -- (4+24/15,0.6);
				\draw[dashed] (1+13*5/15,0) -- (4+26/15,0.4);
				\draw[dashed] (1+14*5/15,0) -- (4+28/15,0.2);

				
			\end{scope}
		\end{tikzpicture}
		
	\end{center}

	\caption{Illustration of the construction of $R_{Int,f}^2$ for a simplex $\Delta = \langle v_1,v_2,v_3 \rangle$ in dimension $2$. The shaded regions do not belong to $R_{Int,f}^2$, there are $4$ points left in $R_{Int,f}^2$.}
\end{figure}
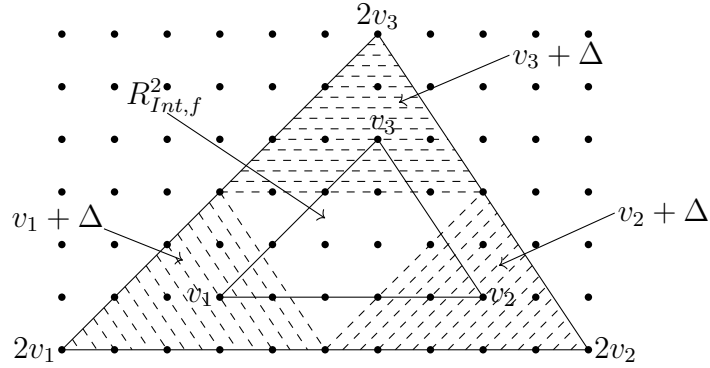

\begin{construction} \label{construction_of_R_Delta_Int_2}
	\normalfont	
	Let $\Delta$ be an $n$-dimensional lattice polytope with given 
	\[  f = \sum\limits_{m \in M \cap \Delta} a_m \cdot x^m. \]
	Given a facet $\Gamma = \Gamma_i \leq \Delta$, where $i \in \{1,...,r\}$, define
	\begin{align*}
		g_{\Gamma}(f) := \sum\limits_{m \in M \cap \Delta} a_m \cdot (\langle n_{i}, m \rangle + b_{i}) \cdot x_0 x^m  \underbrace{=}_{(\textrm{Trafo.})} + b_{i} F_0 + \sum\limits_{j = 1}^n (n_{i})_j F_j \in J_{\Delta,f}^1,
	\end{align*}
	where $n_{i} = ((n_{\Gamma})_1,..., (n_{\Gamma})_n)$ and $b_i = -\Min_{\Delta}(n_i)$. The first representation implies that
	\begin{align*}
		\Supp(g_{\Gamma}(f)) \subset \big( \Delta \cap M \big) \setminus \big( \Gamma \cap M \big).
	\end{align*}
	In case $\Delta$ is an $n$-simplex and
	\[ \Supp(f) = \textrm{Vert}(\Delta)  \]
	then $g_{\Gamma_i}(f) \approx x^{v_i}$ where $v_i$ is the vertex opposite to $\Gamma_i$ ($\approx$ means up multiplication with a nonzero scalar). Conversely the matrix
	\begin{align} \label{matrix_full_rank}
		\begin{pmatrix}
			b_{1} & (n_{\Gamma_1})_1 & ... & (n_{\Gamma_1})_n \\
			\vdots & \vdots \\
			b_{r} & (n_{\Gamma_r})_1 & ... & (n_{\Gamma_r})_n
		\end{pmatrix}	
	\end{align}
	has rank $n+1$ since $n_{\Gamma_1},...,n_{\Gamma_r}$ span $N_{\mathbb{R}}$ and $(b_{1},...,b_{r}) \neq (0,...,0)$. Thus $F_0,...,F_n$ are linear combinations of $g_{\Gamma_1}(f),...,g_{\Gamma_r}(f)$ and we get new generators
	\begin{align} \label{formula_new_generators_Jac_ideal}
		J_{\Delta,f} = (g_{\Gamma_1}(f),...,g_{\Gamma_r}(f)).
	\end{align} 
\end{construction}	

\begin{remark}
	\normalfont
	Working with \textit{Cox coordinate ring} $\mathbb{C}[x_1,...,x_l]$ graded by the divisor class group and given $f \in S_{\beta}$ with $\beta = (b_1,...,b_l)$ under the natural isomorphism
	\begin{align*}
		& S_{\Delta} \rightarrow \bigoplus_{k=0}^{\infty} S_{k \cdot \beta} \\
		& t_0^k \cdot t^m \mapsto \prod\limits_{i} z_i^{k \cdot b_i + \langle m,n_i \rangle}
	\end{align*}
	from (\cite{BatCo94}) the Laurent polynomial $g_{\Gamma_i}(f)$ maps to $z_i \cdot \frac{\partial f}{\partial z_i}$.
\end{remark}

\begin{proposition} \label{proposition_Hodge_comp_Z_f_arbitary_k_dimensions_and_basis}
	Let $\Delta$ be an $n$-dimensional lattice polytope with $l^*(\Delta) > 0$ and a given $f$. Given $\Gamma_1,...,\Gamma_{n+1} \leq \Delta$ with $n_{\Gamma_1},...,n_{\Gamma_{n+1}}$ affine linear independent. Then 
	\begin{align*}
		L^*(k \cdot \Delta) / U_{f,k} = R_{Int,f}^k \quad k=1,...,n+1
	\end{align*}
	where $U_{f,k}$ denotes the vector space over $\mathbb{C}$ spanned by
	\begin{align}
		& g_{\Gamma_i}(f) \cdot x^v \quad i=1,...,n+1, \quad v \in \Int((k-1) \cdot \Delta) \cap M \label{formula_relation_g_Gamma_i_f_the_1}	 \\
		& g_{\Gamma}(f) \cdot x^v \quad \Gamma \leq \Delta, \quad \quad \quad \quad \, \, \, \, v \in \Int((k-1) \cdot \Gamma) \cap M. \label{formula_relation_g_Gamma_i_f_the_2}
	\end{align}
	For $k=2$ these relations are linearly independent.
\end{proposition}

\begin{proof}
	The inclusion
	\[ U_{f,k} \subseteq \Big( J_{\Delta,f}^k \cap L^{*}(k \Delta) \Big).   \]
	is a consequence of the definition of $J_{\Delta,f}$. \\ \\
	\underline{We show:} All relations
	\begin{align} \label{formula_multiplication_in_k_times_interior_of_Delta}
		h \cdot g \in L^*(k \cdot \Delta), \quad \quad \quad h \in J_{\Delta,f}^1, \quad g \in L((k-1) \cdot \Delta)
	\end{align}
	are of type (\ref{formula_relation_g_Gamma_i_f_the_1}) or (\ref{formula_relation_g_Gamma_i_f_the_2}). Finally we exclude (nontrivial) relations
	\begin{align*}
		\sum\limits_{\Gamma} h_{\Gamma} \cdot g_{\Gamma} \in L^*(k \cdot \Delta), \quad h_{\Gamma} \in J_{\Delta,f}^1, \quad g_{\Gamma} \in L((k-1) \cdot \Delta),
	\end{align*}
	\leavevmode \\ \\
	First if $\Supp(g) \subseteq \Int((k-1) \cdot \Delta) \cap M$ then the equation \ref{formula_multiplication_in_k_times_interior_of_Delta} is a linear combination of the relations (\ref{formula_relation_g_Gamma_i_f_the_1}) and without restriction
	\[ \Supp(g) \subset \Bound((k-1) \cdot \Delta) \cap M.  \]
	Write
	\begin{align*}
		h = \sum\limits_{\Gamma} c_{\Gamma} \cdot g_{\Gamma}(f)
	\end{align*}
	and take $F \leq \Delta$ with
	\[ \Supp(g) \cap \Big( (k-1) \cdot F \cap M \Big) \neq \emptyset,   \]
	(this assumption makes sense since $\Supp(g) \neq \emptyset$). Then
	\begin{align*}
		\Supp(h) \cap (F \cap M) = \emptyset
	\end{align*}
	by (\ref{formula_multiplication_in_k_times_interior_of_Delta}). In case $c_{\Gamma} = 0$ except for $\Gamma = F$, then
	\[ \Supp(g) \subset ((k-1) \cdot F)  \]
	and we get a relation of type (\ref{formula_relation_g_Gamma_i_f_the_2}). Restricting $h$ to $F \cap M$ we get
	\begin{align*}
		h_{\vert{F \cap M}} = \sum\limits_{\Gamma} c_{\Gamma} \cdot \big( -b_{\Gamma} F_{0 \vert{\Gamma_1}} + \sum\limits_{i=1}^n (n_{\Gamma})_i F_{i \vert{\Gamma_1}} \big) \in    \langle F_{0 \vert{\Gamma_1}},...,F_{n \vert{\Gamma_1}} \rangle. 
	\end{align*}
	Expanding and restricting to $F$ this means
	\begin{align} \label{formula_relation_between_the_g_Gamma}
		\sum\limits_{\Gamma} c_{\Gamma} \sum\limits_{m \in F \cap M} a_m \big( \langle n_{\Gamma},m \rangle - b_{\Gamma} \big) x_0x^m = 0.
	\end{align}
	The left hand side in (\ref{formula_relation_between_the_g_Gamma}) equals
	\[ \sum\limits_{\Gamma} c_{\Gamma} \cdot \big( -b_{\Gamma} F_{0 \vert{\Gamma_1}} + \sum\limits_{i=1}^n (n_{\Gamma})_i F_{i \vert{\Gamma_1}} \big) \in    \langle F_{0 \vert{\Gamma_1}},...,F_{n \vert{\Gamma_1}} \rangle.  \]
	The $F_{i}$ are algebraically independent over $\mathbb{C}$ (Theorem \ref{theorem_characterization_Delta_regular}), since $f$ is nondegenerate with respect to $\Delta$. Besides $f_{\vert{\Gamma_1}}$ remains nondegenerate with respect to $\Gamma_1$ (by the definition of nondegeneracy) and it follows that there is (up to scaling) only one relation between $F_{0 \vert{\Gamma_1}},...,F_{n \vert{\Gamma_1}}$, the one of the second type. \\ \\
	Now that we know all relations $h \cdot g$ (\ref{formula_multiplication_in_k_times_interior_of_Delta}) we show that these are all relations: For $k=2$ we know all relations by (\cite[Thm.9.8]{Bat93})
	\begin{align*}
		0 \rightarrow \bigoplus\limits_{\Gamma \leq \Delta \textrm{ facet}} L^*(\Gamma) \rightarrow D_f^2 \rightarrow R_{Int,f}^2 \rightarrow 0
	\end{align*}
	and the first map sends $L^*(\Gamma) \ni x^m \mapsto x^m \cdot g_{\Gamma}(f)$. Further $\dim \, D_f^2 = l^*(2 \cdot \Delta) - (n+1) \cdot l^*(\Delta)$ by (\cite[Thm.2.0.11, Rem.2.0.13, Prop.9.0.4]{Bat93}) and we get linear independent relations. For $k \geq 3$ we have an exact sequence
	\begin{align*}
		0 \rightarrow R_{Int,f}^{n+1-k} \rightarrow R_{\Delta,f}^{n+1-k} \rightarrow \bigoplus\limits_{\Gamma} R_{\Gamma,f \vert{\Gamma}}^{n+1-k},
	\end{align*}
	where $\Gamma$ ranges over the facets of $\Delta$. Let us show that this sequence is exact at the middle term: Let $h \in R_{\Delta,f}^{n+1-k}$. If $h \mapsto 0$ then
	\begin{align*}
		h_{\vert{\Gamma}} = \sum\limits_{\tilde{\Gamma}} g_{\Gamma \cap \tilde{\Gamma}} f_{\vert{\Gamma}} \cdot x^{m_{\Gamma,\tilde{\Gamma}}}
	\end{align*}
	for all facets $\Gamma$. Next
	\begin{align*}
		g_{\Gamma \cap \tilde{\Gamma}} f_{\vert{\Gamma}} = g_{\Gamma \cap \tilde{\Gamma}} f_{\vert{\tilde{\Gamma}}}.
	\end{align*}
	Thus restricting $h$ to $\Gamma \cap \tilde{\Gamma}$ and reducing to an equation of linear independent generators we may assume $x^{m_{\Gamma,\tilde{\Gamma}}} = x^{m_{\tilde{\Gamma},\Gamma}}$. Thus 
	\begin{align*}
		h \in J_{\Delta,f}^{n+1-k} \quad \textrm{modulo } L^*((n+1-k) \cdot \Delta)
	\end{align*}
	The term $R_{Int,f}^{n+1-k}$ contributes to the last term $L^*((n+1-k) \cdot \Delta)$. Dualizing the above sequence we obtain
	\begin{align*}
		\bigoplus\limits_{\Gamma} D_{f \vert{\Gamma}}^{k-1} \rightarrow D_f^{k} \rightarrow R_{Int,f}^{k} \rightarrow 0,
	\end{align*}
	where the first map sends $L^*((k-1) \cdot \Gamma) \ni x^m \mapsto x^m \cdot g_{\Gamma}(f)$. The result follows.
	
	
\end{proof}

\begin{example} \label{example_f_support_vertices_of_simplex_jacobian_ring}
	\normalfont
	Let $\Delta$ be a simplex with vertices $v_0,...,v_n$ and assume that $f$ has support on the vertices of $\Delta$. Then as already noted
	\[  g_{\Gamma_i}(f) = x^{v_i}, \] 
	where $\Gamma_i$ denotes the facet opposite to $v_i$. In this case we get a \textit{monomial} basis for example of $R_{Int,f}^2$ by taking the quotient of $L^*(2 \cdot \Delta)$ by the span of 
	\[ x^{v_i + v} \quad i=0,...,n, \quad  v \in \big( \Int(\Delta) \cup \Int(\Gamma_i) \big) \cap M.  \]
\end{example}

\begin{remark} \label{remark_generators_of_D_f_k_all_relations}
	\normalfont
	For $k> 2$ the polynomials in the Proposition will not be linear independent over $\mathbb{C}$ since we have the trivial relations
	\begin{align} \label{formula_relations_between_g_Gamma}
		g_{\Gamma_i}(f) \cdot g_{\Gamma_j}(f) \cdot x^v - g_{\Gamma_j}(f) \cdot g_{\Gamma_i}(f) \cdot x^v
	\end{align}
	for $v \in L^*((k-2) \cdot \Delta)$. For $D_f^k$ we know the \textit{Hilbert-Poincaré series} 
	\begin{align*}
		\dim \, D_f^k = &l^*(k \cdot \Delta) - (n+1) \cdot l^*((k-1) \cdot \Delta) + \binom{n+1}{2} \cdot l^*((k-2) \cdot \Delta) \\
		&- \binom{n+1}{3} \cdot l^*((k-3) \cdot \Delta) + ...
	\end{align*}
	(\cite[Prop.9.0.4]{Bat93}). In fact we have a resolution
	\begin{align*}
		... \rightarrow \bigoplus_{i,j} D_{f \vert{\Gamma_i \cap \Gamma_j}}^{k-2} \rightarrow \bigoplus_{\Gamma} D_{f \vert{\Gamma}}^{k-1} \rightarrow D_f^k \rightarrow R_{Int,f}^k \rightarrow 0
	\end{align*}
	of $R_{Int,f}^k$, which is shown to be exact the same way as above.
	The relations are necessarily of shape \ref{formula_relations_between_g_Gamma} and similarly for more than two indices $i_1,...,i_l$.
	
\end{remark}

\section{Hodge and mixed Hodge structures}

In this section we give some theoretical background on Hodge and mixed Hodge structures.

\begin{definition}
	Let $H$ be a finite dimensional $\mathbb{Q}$-vector space and $H_{\mathbb{C}} := H \otimes \mathbb{C}$. A Hodge structure of weight $k$ on $H$ is a decomposition 
	\[ H_{\mathbb{C}} = \bigoplus\limits_{p+q = k} H^{p,q},   \]
	such that $H^{p,q} = \overline{H^{q,p}}$, where the bar denotes complex conjugation on $H_{\mathbb{C}}$.
\end{definition}

\begin{remark}
	\normalfont
	To put a Hodge structure of weight $k$ on $H$ is equivalent to the existence of a descending filtration
	\[ H_{\mathbb{C}} = F^0 \supset F^1 \supset ... \supset F^{k+1} = 0  \]
	with $H_{\mathbb{C}} = F^p \cap \overline{F^{k+1-p}}$. This correspondence is given by
	\[ F^i = \bigoplus\limits_{p \geq i} H^{p,k-p}, \quad \textrm{and conversely by} \quad   H^{p,q} = F^p \cap \overline{F^q}.  \]
\end{remark}

\begin{definition}
	Let $H$ be a finite dimensional $\mathbb{Q}$-vector space, $H_{\mathbb{C}} := H \otimes \mathbb{C}$. A mixed Hodge structure (for short: MHS) of weight $k$ on $H$ consists of an ascending (weight) filtration
	\[ 0 = W_{k-1} \subset W_k \subset ... \subset W_{2k} = H \]
	and a descending (Hodge) filtration
	\[ H_{\mathbb{C}} = F^0 \supset F^1 \supset ... \supset F^{k+1} = 0  \]
	such that the Hodge filtration induces a Hodge structure of weight $r$ on $W_r/W_{r-1}$.
\end{definition}

We denote the successive quotients $F^{i+1}/F^i$ by $\Gr_F^i$ and $W_r/W_{r-1}$ by $\Gr_W^r$. A Hodge structure of weight $k$ on $H$ defines also a mixed Hodge structure on $H$ by setting $W_k := H$. These abstract definitions are motivated by the following fact:

\begin{theorem} (\cite{Del75}) \\
	For $\Delta$ an $n$-dimensional lattice polytope, $f \in U_{reg}(\Delta)$, the cohomology groups $H^k(Y, \mathbb{Q})$ are equipped with a Hodge structure of weight $k$ and we have 
	\[ H^{p,k-p}(Y):= H^{p,k-p} \cong H^{k-p}(Y, \Omega_{Y}^p). \]
	Besides $H^k(Z_f, \mathbb{Q})$ carry a mixed Hodge structure of weight $k$.
\end{theorem}

If $H = H^k(Y,\mathbb{Q})$ we write $F^iH^k(Y,\mathbb{C})$ and $W_j H^k(Y, \mathbb{Q})$ for $F^i$ and $W_j$.

\begin{definition}
	The dimensions 
	\[ h^{p,k-p}(Y) := \dim \, H^{p,k-p}(Y) \] 
	are called the \textit{Hodge numbers} of $Y$. The spaces $H^{p,k-p}(Y)$ and 
	\[ H^{i,j+k-i}H^k(Z_f, \mathbb{C}) :=  \Gr_F^i \Gr_W^{j+k} H^k(Z_f, \mathbb{C}),\]  
	are called the Hodge components of $Y$ and $Z_f$.
\end{definition}

\begin{definition}
	Let $H, H'$ be $\mathbb{Q}$-vector spaces carrying mixed Hodge structures of weight $k$ and $k+2 \cdot l$ given by filtrations 
	\[ (F^i,W_j) \quad \textrm{and} \quad  ((F')^i,W'_j). \] 
	A $\mathbb{Q}$-linear homomorphism $\phi: H \rightarrow H'$ is called a morphism of mixed Hodge structures of type $(l,l)$ if
	\[ \phi(W_r) \subset W'_{r+2 l}, \quad \phi(F^p) \subset (F')^{p+l}.  \]
\end{definition}

\begin{remark} \label{remark_pullback_morphism_of_MHS}
	\normalfont
	For $Y$ a smooth minimal model of $Z_f$ the inclusion $j: Z_f \rightarrow Y$ induces a pullback homomorphism 
	\[ j^*: H^k(Y, \mathbb{C}) \rightarrow H^k(Z_f, \mathbb{C}). \]
	If we equip $H^k(Y, \mathbb{Q})$ with the trivial weight filtration 
	\[ W_k H^k(Y,\mathbb{Q}) = H^k(Y, \mathbb{Q}), \] 
	then $j^*$ gets a morphism of MHS of type $(0,0)$ (\cite[Ch.7]{Voi02}). This means that the restrictions of cohomology classes of $Y$ to $Z_f$ lie in the minimal weight subspace $W_k H^k(Z_f,\mathbb{C})$. Cohomology classes in $W_{k+j} H^k(Z_f,\mathbb{C})$ for $j>0$ vanish on the comactification $Y$ of $Z_f$.
\end{remark}

\section{The Hodge components of $Z_f$ and $Y_f$}
\label{Section_cohomology_affine_hypersurface_nonsingular_compactification}

\begin{definition}
	We define the primitive cohomology of $Z_f$ as follows:
	\[ PH^{n-1}(Z_{f}, \mathbb{C}) := \coker(H^{n-1}(T, \mathbb{C}) \rightarrow H^{n-1}(Z_{f}, \mathbb{C})). \]
\end{definition}

\begin{remark} \label{remark_Lefschetz_theorems_hypersurfaces_in_tori}
	\normalfont
	There are induced filtrations $F^i$ and $W_j$ on $PH^{n-1}(Z_f, \mathbb{C})$. The torus $T$ is homotopy equivalent to a topological torus, thus
	\[ H^i(T, \mathbb{C}) \cong \wedge^i M_{\mathbb{C}}.  \]
	In particular $h^i(T,\mathbb{C}) = \binom{n}{i}$. By the Lefschetz theorem for hypersurfaces in tori (\cite[Remark 3.10]{DK86}) the restriction homomorphism
	\[ H^{i}(T, \mathbb{C}) \rightarrow H^{i}(Z_f, \mathbb{C})  \]
	is an isomorphism for $i < n-1$ and injective for $i=n-1$.
	
\end{remark}

\begin{remark} \label{remark_hodge_weight_filtration}
	\normalfont
	We have an isomorphism of $\mathbb{C}$ vector spaces (\cite[Thm.6.9, Cor.6.10]{Bat93})
	\[ PH^{n-1}(Z_f, \mathbb{C}) \cong R_{f}.  \]
	This allows us to transport the Hodge and the weight filtration to $R_{f}$ (see \cite[Thm.6.9, Thm.8.2]{Bat93}): The Hogde filtration is given by the reverse grading on $R_{f}$ and the weight filtration on $R_{f}^k$ is induced by the subdivision of $k \cdot \Delta$ into $j$-dimensional faces for $j=0,...,n$. In this thesis we just need the following result:
\end{remark}

\begin{theorem} \label{theorem_batyrev_hodge_components} (\cite[Prop.9.2]{Bat93}) \\
	The Hodge components 
	\[ H^{p,n-1-p}H^{n-1}(Z_f,\mathbb{C}) := Gr_F^p Gr_W^{n-1}H^{n-1}(Z_f,\mathbb{C}) \] 
	are computed as follows:
	\[ H^{p,n-1-p} H^{n-1}(Z_f, \mathbb{C}) \cong R_{Int,f}^{n-p}.  \]
\end{theorem}

\begin{remark}
	\normalfont
	(see \cite{Bat22}, \cite[section 2]{Gie22})\\
	Let $\Delta$ be an $n$-dimensional lattice polytope with $F(\Delta) \neq \emptyset$. The the restrictions of prime toric boundary divisors of $\mathbb{P}_{C(\Delta)}$ and the restrictions of toric strata which are exceptional for $\rho \circ \pi: \mathbb{P} \rightarrow \mathbb{P}_{C(\Delta)}$ compactify $Z_f$ to $Y_f$. We denote the divisorial components of all these divisors by $F_1,...,F_r$ and the sum of the components of codimension $\geq 2$ in $Y_f$ by $S$.
	
\end{remark}	

Let $Z_f^*:= Y_f \setminus D$, where
\[ D:= F_1 + ... + F_r.  \]
There is the following \textit{Gysin exact sequence} (compare also \cite[Proof of Thm.3.7]{DK86}):
\begin{align*}
	... \rightarrow H^{n-2}(Y_f,\mathbb{C}) &\rightarrow H^{n-2}(Z_f^*,\mathbb{C}) \overset{r}{\rightarrow} \bigoplus\limits_{i = 1}^{r+s}  H^{n-3}(F_i, \mathbb{C}) \overset{k_{*}}{\rightarrow} H^{n-1}(Y_f, \mathbb{C}) \\
	&  \overset{j^{*}}{\rightarrow} H^{n-1}(Z_f^*, \mathbb{C}) \overset{r}{\rightarrow} \bigoplus\limits_{i = 1}^{r+s}  H^{n-2}(F_i, \mathbb{C}) \rightarrow ....
\end{align*}

The homomorphism $j^*$ comes from the inclusion $j: Z_{f} \xhookrightarrow{} Y_f$, $r$ is the residue map and $k_{*}$ the so called Gysin map. If we equip $H^{n-1}(Y_f, \mathbb{C})$ with the trivial weight filtration then the map $j^{*}$ is a morphism of MHS of type $(0,0)$ (Remark \ref{remark_pullback_morphism_of_MHS}), the Gysin map is of Hodge type $(1,1)$, whereas $r$ is of Hodge-type $(-1,-1)$ (see \cite[Ch.7]{Voi02}).

\leavevmode
\newline

Restricting to $\Gr_W^{n-1}$ we get
\[ \Gr_W^{n-1} H^{n-2}(Y_f,\mathbb{C}) = 0, \quad \Gr_W^{n-1} H^{n-2}(F_i,\mathbb{C}) = 0.  \]
We arrive at
\begin{align}\label{basic_exact_sequence_cohomology_complement}
	0 &\rightarrow Gr_{W}^{n-1}H^{n-2}(Z_f^*,\mathbb{C}) \overset{r}{\rightarrow} \bigoplus\limits_{i = 1}^r  \Gr_W^{n-3} H^{n-3}(F_i, \mathbb{C}) \overset{k_{*}}{\rightarrow} H^{n-1}(Y_f, \mathbb{C}) \\
	& \overset{j^{*}}{\rightarrow} Gr_{W}^{n-1} H^{n-1}(Z_f^*, \mathbb{C}) \rightarrow 0. \nonumber
\end{align}
Comparing $Z_f^*$ with $Z_f$ we immediately get
\begin{align*}
	& H^{n-1}(Z_f, \mathbb{C}) \cong H^{n-1}(Z_f^*,\mathbb{C}).
\end{align*}



\begin{theorem} \label{theorem_Hodge_comp_min_surf_jac_ring}
	Let $\Delta$ be an $n$-dimensional lattice polytope with $F(\Delta) \neq \emptyset$. Let $Y_f$ be a minimal model of $Z_f$. Then
	\begin{align} 
		H^{p}(Y_f,\Omega_{Y_f}^{n-1-p}) \approx R_{Int,f}^{p+1},
	\end{align}
	where $\approx$ means up to some cohomology classes that arise as restrictions from cohomology classes from strata of $\mathbb{P}$ and which stay constant on the whole family $\mathcal{X} \overset{pr_2}{\rightarrow} B$.
\end{theorem}

\qed

\section{Properties of the period map and its differential} \label{section_period_map_natural_family}

\begin{definition} \label{Def_period_map}
	(\cite[Thm.9.3, Ch.10.1.2, Ch.10.1.3]{Voi02}) \\
	Let $\Delta$ be an $n$-dimensional lattice polytope with $l^*(\Delta) > 0$, let $f \in B$ and assume $Y_f$ is smooth. The period map $\mathcal{P}_{B,f}$ for the $n-1$-th cohomology is defined by
	\begin{align*}
		\mathcal{P}_{B,f} = \bigoplus\limits_{k=1}^{n-1} \mathcal{P}_{B,f}^k: &B \rightarrow \Gamma \setminus D \\
		& f' \mapsto (F^{n-1}H^{n-1}(Y_{f'},\mathbb{C}),...,F^1 H^{n-1}(Y_{f'},\mathbb{C}))
	\end{align*}
	where
	\[ H^{n-k,k-1}(Y_{f'}) \subset H^{n-1}(Y_{f'}, \mathbb{C}) \cong H^{n-1}(Y_f,\mathbb{C})  \]
	Here
	\begin{align*}
		& D: \textrm{ period domain  (a quasiprojective variety) } \\
		& \Gamma: \textrm{The monodromy group } \pi_1(B,f)
	\end{align*}
	We refer to (\cite[Ch.10]{Voi02}, \cite[Ch.3]{Voi03}) for details.
\end{definition}

\begin{construction} (Result of Griffiths) (\cite[Thm.10.21]{Voi02}) \label{construction_result_of_Griffiths}
	\normalfont \\
	The differential $d \mathcal{P}_{B,f}$ of $\mathcal{P}_{B,f}$ fits into a diagram
	\begin{equation}
		\begin{tikzcd} \label{Commutative_diagram_Kodaira_Spencer_map_dP}
			T_{B,f}  \arrow[swap]{dr}{d \mathcal{P}_{B,f} } \arrow{r}{\kappa_{f}} & H^{1}(Y_f,T_{Y_f})  \arrow{d}{\Phi_f} \\
			& \bigoplus\limits_{k=0}^{n-2} \Hom(H^{k}(Y_f, \Omega_{Y_f}^{n-1-k}), H^{k+1}(Y_f, \Omega_{Y_f}^{n-2-k}))
		\end{tikzcd}
	\end{equation}
	since $Y_f$ is smooth. $\Phi_{f}$ is the homomorphism between cohomology groups induced by cup product and the contraction
	\[ T_{Y_f} \times \Omega_{Y_f}^{n-1-k} \rightarrow \Omega_{Y_f}^{n-2-k}.  \]
	This diagram is important since it connects the Hodge-theoretic homomorphism $d \mathcal{P}_{B,f}$ with the Kodaira-Spencer map $\kappa_f$. 
\end{construction}

\begin{remark}
	\normalfont
	Starting with a smooth proper deformation $\mathcal{Y} \rightarrow S$ of $Y_f$ with $S$ smooth, we define a Kodaira Spencer map $\kappa_{S,f}$ and a period map $\mathcal{P}_{S,f}$ just as in the Definitions above. The result of Griffiths remains valid: $d \mathcal{P}_{S,f}$ factors through $\kappa_{S,f}$ and $\Phi_f$, that is $\Phi_f$ is universal.	
\end{remark}	

\begin{definition}
	The infinitesimal Torelli Theorem (short: ITT) for $Y_f$ asks if $\Phi_f$ is injective. The infinitesimal Torelli Theorem for $Y_f$ in $\mathcal{X} \overset{pr_2}{\rightarrow} B$ asks if $\Phi_{f \vert{Im \, \kappa_f}}$ is injective. The infinitesimal Torelli Theorem for $\mathcal{X} \overset{pr_2}{\rightarrow} B$ asks if $\Phi_{f \vert{Im \, \kappa_f}}$ is injective for $f \in B$. 
\end{definition}

If $\kappa_f$ is surjective of course the first and the third definition coincide. Given the identifications (\ref{theorem_batyrev_hodge_components}) choose a reference point $f \in B$ and define a \textit{mixed period map} $\phi_f$ from $B$ into a \textit{mixed period domain} $D_{mix}$ (which has to be defined) by
\[ f' \overset{\phi_f}{\mapsto}  (R_{Int,f'}^1,..., R_{Int,f'}^{n}).  \]
Apparently $\phi_f$ just depends on the affine part $Z_f$ and not on the particular compactification $Y_f$. 
\leavevmode
\\

\begin{lemma} \label{lemma_diff_period_map_induces_hom_on_R_Int,f}
	$d \phi_f^k$ has image in $\Hom(R_{Int,f}^k, R_{Int,f}^{k+1})$ and $d \mathcal{P}_{B,f}^k$ factors as follows
	\begin{equation}
		\begin{tikzcd} \label{commutative_diagram_dP_B_mixed_period_map}
			T_{B,f}  \arrow[swap]{dr}{d \mathcal{P}_{B,f}^k } \arrow{r}{d \phi_f^k} & \Hom(R_{Int,f}^k, R_{Int,f}^{k+1})  \arrow{d} \\
			& \Hom(H^{k-1}(Y_f, \Omega_{Y_f}^{n-k}), H^k(Y_f, \Omega_{Y_f}^{n-k+1}))
		\end{tikzcd} 
	\end{equation}
	where the vertical map denotes the inclusion.
\end{lemma}
\qed

\section{Smooth and stable points of $\mathcal{M}(\Delta)$} \label{section_smooth_and_stable_points}

Given an $n$-dimensional lattice polytope $\Delta$ and some polynomial $f$, the torus $T$ acts on $U_{reg}(\Delta)$
\[ (t_1,...,t_n).f(x_1,...,x_n) = f(t_1x_1,...,t_n x_n), \quad (t_1,...,t_n) \in T.  \]
$T \cdot f$ denotes the orbit of $f$ under $T$. \\ \\
Let $\mathcal{M}(\Delta) := B/T$ be the quotient of $B$ by $T$. We omit equivalence classes and write $f \in \mathcal{M}(\Delta)$.

\begin{definition} \label{definition_stable_points} (\cite[Def.1.7]{MuFo82}) \\
	Let $v \in \mathbb{C}^{l(\Delta) +1}$ and $x = [v]$. The point $x$ is  stable if the orbit $T \cdot v$ is closed and of dimension $\dim \, \Delta = n$. The second condition is equivalent to the condition
	\[ \# \Stab_T(v) < \infty, \]
	where 
	\[  \Stab_T(v) := \{ t \in T| \, t.v = v \} \]
	denotes the stabilizer of $v$ w.r.t. the action of $T$.
\end{definition}

Let $r:= l(\Delta)-1$. The set $(\mathbb{P}^r)^s$ of stable points of $\mathbb{P}^r$ is Zariski open in $\mathbb{P}^r$ (\cite[§4]{MuFo82}), but might be empty. Let 
\[ \Delta_r := \langle e_m \, | \quad  m \in M \cap \Delta \rangle \] 
denote the $r$-dimensional standard simplex embedded into an affine hyperplane in $\mathbb{R}^{l(\Delta)}$. Then there is a map $\pi: \Delta_r \rightarrow \Delta$ given by
\[ \sum\limits_{m \in M \cap \Delta} \lambda_m \cdot e_m \mapsto \sum\limits_{m \in M \cap \Delta} \lambda_m \cdot m \quad \textrm{for } \sum\limits_{m \in M \cap \Delta} \lambda_m = 1. \]
To $\Delta_r$ is associated the toric variety $\mathbb{P}^r$ and given 
\[ a:= (a_m)_{m \in M \cap \Delta} \in \mathbb{P}^r \] 
there is a natural $(\mathbb{C}^*)^r$ orbit through a of some dimension $k \in \{0,...,r\}$. We denote the $k$-dimensional face of $\Delta_r$ corresponding to this orbit by $\Gamma(a)$.

\begin{proposition} (\cite[Prop.3.5]{KSZ91}) \\
	$(\mathbb{P}^r)^s \neq \emptyset$ if and only if 
	\[  (0,...,0) \in \Int(\Delta) \cap M. \] 
	More precisely $a \in (\mathbb{P}^r)^s$ is stable if and only if $\pi(\Gamma(a))$ has full dimension $n$ and contains $(0,...,0)$ in its interior.
\end{proposition}

\begin{corollar} \label{corollary_f_nondegenerate_then_f_stable}
	If $(0,...,0) \in \Int(\Delta) \cap M$ and $a:= (a_m)_{m \in M \cap \Delta} \in \mathbb{R}^{l(\Delta)}$ is such that $f$  lies in $ U_{reg}(\Delta)$, then $a$ is stable.
\end{corollar}

\begin{proof}
	For $m$ a vertex of $\Delta$ we have $a_m \neq 0$ and thus $\pi(\Gamma(a)) \supset \Int(\Delta)$, in particular $(0,...,0) \in \pi(\Gamma(a))$ and $\pi(\Gamma(a))$ is full-dimensional.
\end{proof}


\begin{corollar} \label{corollary_stable_nondeg_point_automatically_smooth}
	If $(0,...,0) \in \Int(\Delta)$ then the quotient $\mathcal{M}(\Delta)$ is smooth.
\end{corollar}

\begin{proof}
	By definition $L(\Delta) \setminus U_{reg}(\Delta) = \{E_A = 0\}$, where $A := \Delta \cap M$. In effect $U_{reg}(\Delta)$ is affine.
	\\ \\
	The orbit $T \cdot f$ is closed and $n$-dimensional since all $f \in \mathcal{M}(\Delta)$ are stable w.r.t. the action of $T$ on $\mathbb{P}^{l(\Delta) -1}$. Applying Luna's slice Theorem (\cite[Prop.5.7]{Dre12}) to the affine variety $U_{reg}(\Delta)$ and the projection
	\[ U_{reg}(\Delta) \rightarrow \mathcal{M}(\Delta)  \]
	gives us that $\mathcal{M}(\Delta)$ is smooth at $f$ if $\Stab_{T \times \mathbb{C}^*}(f)$ contains just the neutral element $(1,...,1)$. \\ \\
	By definition
	\[ \Stab_T(f) = \{ (t,1) = (t_1,...,t_n,1) \in T \times \mathbb{C}^*  \, | \quad t^m  a_m = a_m  \quad \forall m \, \textit{ with } a_m \neq 0  \}.  \]
	Thus
	\[ \Stab_T(f) \subset  \{ t \in T \, | \, t^{v_i} = 1 \quad \forall \textrm{ vertices } v_1,...,v_k   \} \times \{1\}^. \]
	\leavevmode
	\\
	Consider $n$ vertices $v_1,...,v_n$ which span $M_{\mathbb{R}}$ and apply an unimodular transformation $U: M \rightarrow M$ such that 
	\[ U(v_1),...,U(v_{n-1}) \in \{ (m_1,...,m_n) \in M| \, m_n = 0 \}. \] 
	Replace $v_i$ by $U(v_i)$ and note that for
	\[ (t_1,...,t_n) \in \Stab_T(f) \] 
	the entry $t_n$ is uniquely determined by $t_1,...,t_{n-1}$ and the relation
	$1 = t^{v_n}$. Further if $(t_1,...,t_{n-1}) = (1,...,1)$ then also $t_n = 1$. By this we have reduced the assertion $\Stab_T(f) = \{(1,...,1)\}$ to the lower-dimensional problem that the only solution of
	\[ 1 = t^{v_1} = ... = t^{v_{n-1}}  \]
	is $t= (1,...,1)$. We continue inductively.
	
\end{proof}

\begin{remark} \label{remark_tangent_space_M(Delta)}
	\normalfont
	Given $(0,...,0) \in \Int(\Delta) \cap M$ and $f \in \mathcal{M}(\Delta)$ then $f$ is smooth and stable. Let $T \rightarrow B$ denote the map
	\[ (t_1,...,t_n) \mapsto f(t_1x_1,...,t_n x_n).  \]
	Then we get
	\begin{align*}
		T_{\mathcal{M}(\Delta),f} &\cong T_{B,f}/T_{T \cdot f, f} \\
		&\cong L(\Delta)/ \langle f,  x_1 \frac{\partial f}{\partial x_1},..., x_n \frac{\partial f}{\partial x_n} \rangle \cong R_{f}^1
	\end{align*} 
	where the first isomorphism follows from (\cite[Cor.11.3]{Bat93}).
\end{remark}

\begin{remark}
	\normalfont
	Given $f \in \mathcal{M}(\Delta)$ the differential $d \phi_f$ and $\kappa_f$ factor through $T_{\mathcal{M}(\Delta),f}$. The differential $d \phi_f$ is simply induced by the addition of lattice points (see \cite[Prop.11.8]{Bat93})
	\begin{equation}
		\begin{tikzcd}
			& L(\Delta) \arrow{r}{+} \arrow{d} &  \Hom(L^*(k \cdot \Delta), L^*((k+1) \cdot \Delta)) \arrow{d} \\
			& R_f^1 \arrow{r}{d \phi_f^k} & \Hom(R_{Int,f}^k, R_{Int,f}^{k+1}))
		\end{tikzcd}
	\end{equation}
\end{remark}

\section{Explicit description of $\ker(d \phi_{f})$}

The following elementary Remark follows also from diagram \ref{Commutative_diagram_Kodaira_Spencer_map_dP}.

\begin{remark} \label{remark_kernel_kappa_f_subset_ker_dP_B_f}
	(Elementary proof that $\ker(\kappa_f) \subset \ker(d \phi_{f}^k$)) \\
	\normalfont
	We use the notation from (\cite{Gie22}): Let $\alpha \in R(N,\Sigma_{\Delta})$. Then there is $\Gamma_{-\alpha} \leq \Delta$ such that if $v \in \Int(k \cdot \Delta) \cap M$ then by definition of roots 
	\[ v-\alpha \in \Int(k \cdot \Delta) \cap M \quad \textrm{or} \quad v-\alpha \in \Int(k \cdot \Gamma_{-\alpha}) \cap M.  \]
	We get 
	\[ g_{\Gamma_{-\alpha}}(f) \cdot x^{-\alpha} = w_{-\alpha}(f).  \]
	and
	\[ w_{-\alpha}(f) \cdot x^v = g_{\Gamma_{-\alpha}}(f) \cdot x^{v-\alpha} \in L^*((k+1) \cdot \Delta).  \]
	Since $g_{\Gamma_{-\alpha}}(f) \cdot x^{v-\alpha} \in J_{\Delta,f}^{k+1}$ we verified explicitly that $w_{-\alpha}(f) \in \ker(d \phi_{f}^k)$ for all $k=1,...,n-1$. 
\end{remark}


By definition 
	\begin{align} \label{general_formula_for_ker_DP_MDelta_not_handable}
		\ker(d \phi_{f}^k) = \{ h \in L(\Delta) \mid \, \forall v \in \Int(k \cdot \Delta) \cap M: \, h \cdot x^v \in J_{\Delta,f}^{k+1} \}.
	\end{align}

\begin{theorem} \label{Theorem_formula_ker_ITT}
	Let $\Delta$ be an $n$-dimensional lattice polytope 
	with $(0,...,0) \in \Int(\Delta) \cap M$. If 
	$f \in L(\Delta)$ is nondegenerate
	\begin{align}\label{general_formula_for_ker_DP_MDelta_handable}
		\ker(d \phi_{f}^1) =&  \Big\langle g_{\Gamma}(f) \cdot x^w \in R_f^1 \mid \Gamma \leq \Delta \textrm{ a facet}, \nonumber  \\
		& w + v \in \big( \Int(\Delta) \cup \Int(\Gamma) \big) \cap M, \\
		& \forall \, v \in \Int(\Delta) \cap M \Big\rangle. \nonumber
	\end{align}
\end{theorem}

The inclusion $\supseteq$ is obvious. The main problem with the opposite inclusion is that it might happen that 
\begin{align} \label{formula_difficulty_kernel_linear_comb_lDelta}
	\sum\limits_{\Gamma} h_{\Gamma} \cdot g_{\Gamma}(f) \in L(\Delta).
\end{align}
but $h_{\Gamma} \cdot g_{\Gamma}(f) \notin L(\Delta)$ for at least two $\Gamma$'s.

\begin{remark}
	\normalfont
	Given $\Delta$ with $(0,..,0) \in \Int(\Delta) \cap M$ and $g_{\Gamma}(f) \cdot x^w \in \ker(d \phi_f^1)$ assume that we replace $M$ by some $p \in M$ such that the new zero point $-p$ is again an interior lattice point of the new polytop $\tilde{\Delta}$. Then $g_{\Gamma}(f) \cdot x^w$ gets well transformed to $g_{\tilde{\Gamma}}(f) \cdot x^w$ for $\tilde{\Gamma} \subset \tilde{\Delta}$, where $g_{\tilde{\Gamma}}(f) \cdot x^{-p} = g_{\Gamma}(f)$, with equivalent relation \ref{general_formula_for_ker_DP_MDelta_handable}. This is basically because $-p \in \Int(\Delta) \cap M$ and due to the relation \ref{general_formula_for_ker_DP_MDelta_handable} again.	
\end{remark}


	Set $R_{\Gamma} := (\Int(\Delta) \cup \Int(\Gamma)) \cap M$ and write
	\begin{align*}
		h = h \cdot x^{(1,0,...,0)} \cdot x_0^{-1} = \sum h_{\Gamma} \cdot g_{\Gamma}(f) \cdot x_0^{-1} \in \ker(d \phi_f^1),
	\end{align*}
	where $\Supp \, h_{\Gamma} \subset R_{\Gamma}$ by \ref{general_formula_for_ker_DP_MDelta_not_handable}. Consider the following three points

\begin{enumerate}
	\item $h_{\Gamma} \cdot g_{\Gamma}(f) \in \Hom(R_{Int,f}^k, R_{Int,f}^{k+1}), \quad k=1,n-1$
	\item $m+v \in R_{\Gamma}$ for all $m \in \Supp(h_{\Gamma})$ and $v \in \Int(\Delta) \cap M$.
	\item $h_{\Gamma} \cdot g_{\Gamma}(f) \in L(\Delta)$.
\end{enumerate}

Below we show $2.$, which is a strengthening of $1.$ for $k=1$. The case $k = n-1$ in $1.$ follows by duality:
\begin{align*}
	\Hom(R_{Int,f}^1,R_{Int,f}^2) \cong \Hom(R_{Int,f}^{n-1},R_{Int,f}^n).
\end{align*}
Then we show that $1$. for $k=n-1$ is equivalent to $3.$. Finally we may apply Lemma \ref{lemma_proof_ITT_just_one_supp} to $3.$ to finish the proof of theorem \ref{Theorem_formula_ker_ITT}.
\leavevmode
\\ \\
To $2.$: Definitely we always have $h_{\Gamma} \cdot x^v \in L^*(2 \cdot \Delta)$ for $v \in \Int(\Delta) \cap M$ and thus we may represent 
\begin{align} \label{equation_D_f_3_relation_h_Gamma_Gamma_v}
	D_f^3 \ni 0 \equiv \sum h_{\Gamma} \cdot x^v \cdot g_{\Gamma}(f) = \sum h_{\Gamma,v} \cdot g_{\Gamma}(f),
\end{align}
as relation between the generators of zero in $D_f^3$, where $\Supp(h_{\Gamma,v}) \subset R_{\Gamma}$. By Remark \ref{remark_generators_of_D_f_k_all_relations} either 
\begin{align} \label{equation_h_Gamma_beginning_form}
	h_{\Gamma,v} = h_{\Gamma} \cdot x^v,
\end{align}
which implies $\Supp \, h_{\Gamma} \cdot x^v = \Supp(h_{\Gamma,v}) \subset R_{\Gamma}$ and $v$ causes no problems, or else there is a facet $\Gamma_2 \neq \Gamma_1:= \Gamma$, $v_1:=v$, such that
\begin{align}
	& h_{\Gamma_1,v_1} - h_{\Gamma_1} \cdot x^{v_1} = c_2 \cdot g_{\Gamma_2}(f) \cdot x^{v_2} \label{equation_relation_in_D_f_k_first} \\
	& h_{\Gamma_2,v_1} - h_{\Gamma_2} \cdot x^{v_1} = c_2 \cdot g_{\Gamma_1}(f) \cdot x^{v_2}, \label{equation_relation_in_D_f_k_second}
\end{align}
where $c_2 \neq 0$, $v_2 \in \Int(\Delta) \cap M$ and $\Supp \, h_{\Gamma_i,v_1} \subset R_{\Gamma_i}$. These relations are the only possible relations: For if we if we insert a new term
\begin{align*}
	&g_{\Gamma_2}(f) \cdot x^{v_2} \cdot g_{\Gamma_1}(f) - g_{\Gamma_1}(f) \cdot x^{v_2} \cdot g_{\Gamma_2}(f) + g_{\Gamma_3}(f) \cdot x^{v_3} \cdot g_{\Gamma_2}(f) - \\
	&g_{\Gamma_2}(f) \cdot x^{v_3} \cdot g_{\Gamma_3}(f),
\end{align*}
which cancels, then $g_{\Gamma_1}(f) \cdot x^{v_2} = g_{\Gamma_3}(f) \cdot x^{v_3}$ and we achieve a desired relation again. We skip the nonzero constant $c_2$ from the notation.
\leavevmode 
\\ \\
Obviously $g_{\Gamma_2}(f) \cdot x^{v_2} \notin L(\Delta)$ (and similarly for the second equation) for else $h_{\Gamma_1} \cdot x^{v_1} \in L(\Delta)$, $c_2 = 0$ and $h_{\Gamma_1,v_1} = h_{\Gamma_1} \cdot x^{v_1}$ by Proposition \ref{proposition_Hodge_comp_Z_f_arbitary_k_dimensions_and_basis} (the relations are linear independent for $k=1$). In particular $v_2 \neq (0,...,0)$.
Next note that 
\begin{align} 
	& \langle v_2,n_F \rangle > \Min \big( 0, \langle v_1,n_F \rangle \big) \quad \forall F \neq \Gamma_1,\Gamma_2, \quad \label{equation_finitely_many_steps_scalar_product_gets_positiver} \\
	& \langle v_2,n_{\Gamma_i} \rangle \geq \langle \Min \big( 0, v_1,n_{\Gamma_i} \rangle \big)  \quad  \Gamma_i = \Gamma_1,\Gamma_2. \label{equation_finitely_many_steps_scalar_product_gets_positiver_second}
\end{align}
by \ref{equation_relation_in_D_f_k_first} and \ref{equation_relation_in_D_f_k_second}. We write $v_2 \geq_F v_1$ for these relations. 
Write $v_2 >_F v_1$ if $v_2 \geq_F v_1$ and if there is a facet $F$ with $0 > \langle v_1,n_F \rangle < \langle v_2,n_F \rangle$. \leavevmode \\ \\
Assume \ref{equation_relation_in_D_f_k_first} and \ref{equation_relation_in_D_f_k_second} are true for some $v_1 \in \Int(\Delta) \cap M$ (with $c_2 \neq 0$): We may choose $v_1$ in such a way such that given $w \in \Int(\Delta) \cap M$ with $w >_F v_1$ then
\begin{align*}
	h_{\Gamma_1,w} = h_{\Gamma_1} \cdot x^{w}, \quad h_{\Gamma_2,w} = h_{\Gamma_2} \cdot x^{w}.
\end{align*}

\leavevmode \\
Besides we may choose $v_1$ such that $v_2 >_F v_1$: Obviously $v_2 \geq_F v_1$ and if there is a facet $F \neq \Gamma_1,\Gamma_2$ with $\langle v_1,n_F \rangle < 0$ then $v_2 >_F v_1$ by \ref{equation_finitely_many_steps_scalar_product_gets_positiver}. Else assume for example
\begin{align*}
	0 > \langle v_1,n_{\Gamma_1} \rangle = \langle v_2,n_{\Gamma_1} \rangle.
\end{align*}
Then $\langle v_2,n_F \rangle > \langle v_1,n_F \rangle$ for all facets $F \neq \Gamma_2$ with $F \cap \Gamma_1 \neq \emptyset$ by \ref{equation_relation_in_D_f_k_first}. Up to replacing $v_1$ by $v_2$ with the bounded condition $v_1,v_2 \in \Int(\Delta) \cap M$ after finitely many replacements we get a vector $v_1$ with $\langle v_1,n_F \rangle$ maximal among those vectors $w \in \Int(\Delta) \cap M$ with $\langle w,n_{\Gamma_1} \rangle = \langle v_1,n_{\Gamma_1} \rangle$, forcing $v_2 >_F v_1$ henceforth. Then by assumption
\begin{align}
	& L(\Delta) \ni h_{\Gamma_1,v_2} = h_{\Gamma_1} \cdot x^{v_2} = \big( h_{\Gamma_1,v_1} - g_{\Gamma_2}(f) \cdot x^{v_2} \big) \cdot x^{v_2-v_1} 
	\label{equation_resolve_v_2_v2-v1_first}
	\\
	& L(\Delta) \ni h_{\Gamma_2,v_2} = h_{\Gamma_2} \cdot x^{v_2} = \big( h_{\Gamma_2,v_1} - g_{\Gamma_1}(f) \cdot x^{v_2} \big) \cdot x^{v_2-v_1}.
	\label{equation_resolve_v_2_v2-v1_second}
\end{align}
Thus if $\langle v_2,n_F \rangle < 0$ then $\langle 2 \cdot v_2-v_1,n_F \rangle \geq 0$, that is $2 \cdot v_2 \in \Int(\Delta) \cap M$ and $2 \cdot v_2 \geq_F v_1$. 

\leavevmode \\
By the same reasoning as above but arguing with equation \ref{equation_resolve_v_2_v2-v1_first}, \ref{equation_resolve_v_2_v2-v1_second}  instead of equation \ref{equation_relation_in_D_f_k_first}, \ref{equation_relation_in_D_f_k_second} we get $2 \cdot v_2 >_F v_1$. Thus
\begin{align*}
	& L(\Delta) \ni h_{\Gamma_1,2v_2} = h_{\Gamma_1} \cdot x^{2v_2} = \big( h_{\Gamma_1,v_1} - g_{\Gamma_2}(f) \cdot x^{v_2} \big) \cdot x^{2v_2-v_1} 
	\\
	& L(\Delta) \ni h_{\Gamma_2,2v_2} = h_{\Gamma_2} \cdot x^{2v_2} = \big( h_{\Gamma_2,v_1} - g_{\Gamma_1}(f) \cdot x^{v_2} \big) \cdot x^{2v_2-v_1}.
\end{align*}
We conclude directly $3 \cdot v_2 \in \Int(\Delta) \cap M$ and $3 \cdot v_2 \geq_F v_1$ and by the same argument as above $3 \cdot v_2 >_F v_1$. In this way we would get an infinite series $k \cdot v_2 >_F v_1$ which is absurd for $v_2 \neq (0,...,0)$.

\leavevmode 
\\
Thus \ref{equation_h_Gamma_beginning_form} is true for all $\Gamma$ and all $v \in \Int(\Delta) \cap M$ and we proved $2.$.

\begin{lemma} \label{lemma_proof_ITT_just_one_supp}
	Assume that
	\[ h_{\Gamma} \cdot g_{\Gamma}(f) \in L(\Delta)  \]
	for some facet $\Gamma$. Then 
	\begin{align} \label{formula_support_h_Gamma_other_normal_vectors_3}
		x^m \cdot g_{\Gamma}(f) \in L(\Delta) \quad \forall \, m \in \Supp(h_{\Gamma})
	\end{align}
\end{lemma}

\begin{proof} 
	Assume to the contrary that $x^m \cdot g_{\Gamma}(f) \notin L(\Delta)$, say that this polynomial jumps out of a facet $\Gamma_1$. Then choose $n \in N$ suitable and find a vertex $v \in Vert(\Gamma_1)$ such that $m+v \notin \Delta \cap M$ and this vector remains left in $\Supp(h_{\Gamma} \cdot g_{\Gamma}(f))$.
\end{proof}


	The point 1. is true for $k=n-1$. This means that
	\begin{align} \label{equation_condition_R_Int_f_n}
		h_{\Gamma} \cdot x^{w} \cdot g_{\Gamma}(f) \in L^*(n \cdot \Delta) \quad \forall w \in L^*((n-1) \cdot \Delta).
	\end{align}

	\begin{lemma}
		We have $h_{\Gamma} \cdot g_{\Gamma}(f) \in L(\Delta)$.
	\end{lemma}
	
\begin{proof}
	Assume $\langle m, n_{F} \rangle < 0$ for some facet $F \neq \Gamma$. Take the $n$-dimensional convex lattice polytope 
	\begin{align*}
		P_{F} := \Big\langle (0,...,0), F \Big\rangle
	\end{align*}
	and remember that $(0,...,0) \in \Int(\Delta) \cap M$. By combinatorics of lattice polytopes $(n-1) \cdot P_{F}$ is normal (\cite[Thm. 2.2.12]{CLS11}), thus for $k \in \mathbb{N}$
	\begin{align*}
		(k \cdot (n-1) \cdot P_{F}) \cap M = \underbrace{(n-1) \cdot P_{F} \cap M + ... + (n-1) \cdot P_{F} \cap M}_{k \textrm{ summands}}
	\end{align*}
	Choosing $k \gg 0$ there is $x \in k \cdot (n-1) \cdot P_{F} \cap M$ with$\langle x,n_{F} \rangle = k \cdot (n-1) \cdot \Min_{\Delta}(n_{F}) + 1$. By normality 
	\begin{align*}
		x = m_1 + ... + m_k, \quad \textrm{ where } m_i \in (n-1) \cdot P_{F} \cap M
	\end{align*}
	and up to permutation $m_i \in (n-1) \cdot F \cap M$ for $i=1,...,k-1$ and we get an element $m_k \in \Int((n-1) \cdot \Delta) \cap M$ with $\langle m_k,n_{F} \rangle = (n-1) \cdot n_{F} +1$. \\ 
	We obtain $\langle m + m_k, n_{F} \rangle \leq (n-1) \cdot \Min_{\Delta}(n_{F})$, thus $x^{m+m_k} \cdot g_{\Gamma}(f) \notin L^*(n \cdot \Delta)$, conradicting \ref{equation_condition_R_Int_f_n} for $w:= m_k$. Thus $\langle m,n_F \rangle \geq 0$ for all facets $F \neq \Gamma$. We get 
	\begin{align} \label{formula_zwischenformula_hGamma_in_LDelta}
		h_{\Gamma} \cdot g_{\Gamma}(f) \in L(\Delta)
	\end{align}
	since $h \in L(\Delta)$ and since the term \ref{formula_zwischenformula_hGamma_in_LDelta} jumps at most out of the facet $\Gamma$ (different terms do not compensate each other).
\end{proof}

The theorem has strong consequences concerning the ITT: Suppose given an element of the kernel as in \ref{general_formula_for_ker_DP_MDelta_handable}. Let 
\[ H_{n,l} := \{ x \in M_{\mathbb{R}} \mid \langle x,n \rangle = l \} \quad n \in N, \quad l \in \mathbb{Z}.  \]
\begin{corollar} \label{lemma_def_l_i_Supp_h_Gamma_i}
	Given the assumptions of theorem \ref{Theorem_formula_ker_ITT} we have
	\begin{align} \label{support_h_Gamma_l_Gamma}
		& \Supp(h_{\Gamma}) \subset \Cone \big(\Int(\Gamma) \cap M \big) \cap H_{n_{\Gamma},-l_{\Gamma}} \cap M,
	\end{align}
	where $l_{\Gamma}$ denotes the smallest natural number $\geq 1$ with 
	\[ \emptyset \neq H_{n_{\Gamma}, -l_{\Gamma}} \cap R_{\Gamma} \]
	(if $l_{\Gamma}$ does not exist, then $\Supp(h_{\Gamma}) = \emptyset$).
\end{corollar}

\begin{proof}
	If 
	\[ (0,0,0) \neq m \in \Supp(h_{\Gamma})  \]
	then some multiple $r \cdot m$, $r \in \mathbb{N}_{\geq 1}$, lies in $\Int(\Gamma) \cap M$ by \ref{general_formula_for_ker_DP_MDelta_handable}. Thus $m \in \Cone(\Int(\Gamma) \cap M)$. If there were $m' \in \Int(\Delta) \cap M$ with
	\[ 0 > \langle m', n_{\Gamma} \rangle > \langle m,n_{\Gamma} \rangle  \]
	then
	\begin{align*}
		& (r - 1) \cdot m + m' \in \Int(k \cdot \Delta) \cap M \\
		& \Rightarrow \langle  r \cdot m + m', n_{\Gamma} \rangle = \Big( \underbrace{(r-1) \cdot m + m'}_{\in \Int(\Delta) \cap M} \Big) + \underbrace{m}_{\in \Supp(h_{\Gamma})}  > \Min_{\Delta}(n_{\Gamma})
	\end{align*}
	a contradiction.	
\end{proof}

\begin{example}
	\normalfont	
	Take $\mathbb{P} = \mathbb{P}^n$, where $n \geq 3$ and $Y_f$ a generic nondegenerate (not exactly the same as smooth) hypersurface in $\mathbb{P}^n$ of degree $d \geq n+2$. Then the assumptions of theorem (\ref{Theorem_formula_ker_ITT}) are met and the Kodaira-Spencer map is surjective (see \cite[Lemma 6.15]{Voi03}). Thus we get a different proof of the ITT for such projective hypersurfaces.
\end{example}

\begin{corollar} \label{corollary_inf_Tor_thm}
	Assume $\Int(\Delta) \cap M \nsubseteq \{plane\}$ and $f$ is nondegenerate. Then $\ker(d \phi_f) = \ker(\kappa_f)$. Thus
	\[  \ker(\Phi_{f \vert{Im \, \kappa_f}}) = \{ 0  \},  \]
	that is $Y_f$ fulfills the ITT for $\mathcal{X} \rightarrow B$. 
\end{corollar}

\begin{proof}
	We show that $\ker(d \phi_f^1) = \ker(\kappa_f)$. The corollary follows from this assertion since $\ker(\kappa_f) \subset \ker(d \phi_f^k)$ for $k=1,...,n-1$ by Remark \ref{remark_kernel_kappa_f_subset_ker_dP_B_f} and thus
	\begin{align*}
		\ker(d \phi_f) = \bigcap\limits_{k} \ker(d \phi_f^k) = \ker(\kappa_f).
	\end{align*}
	Let $f \in B$ and $g \in \ker(d \phi_{f}^1)$. By theorem \ref{Theorem_formula_ker_ITT} we may assume that 
	\[ g = g_{\Gamma}(f) \cdot x^w,  \]
	where
	\[ w \in R_{\Gamma} \cap M  \]
	for some $\Gamma \leq \Delta$ and $\langle w,n_{\Gamma} \rangle = -k \leq 0$, $\langle w, n_{\Gamma'} \rangle \geq 0$ for $\Gamma' \neq \Gamma$. There are two simple cases	\begin{itemize}
		\item $k = 0$:   Then $w = (0,0,0)$ and $g = g_{\Gamma}(f) \equiv 0 \in R_{f}^1$.
		\item $k = -1$:   Then $w \in -R(N,\Sigma_{\Delta})$ and $g \in \ker(\kappa_{f})$.
	\end{itemize}
	Else $k \leq -2$ and there is no lattice point $w' \in \Int(\Delta) \cap M$ with
	\[ 0 > \langle w', n_{\Gamma} \rangle > \langle w, n_{\Gamma} \rangle.  \]
	Choose $x:= (0,0,0)$, and $y,z \in \Int(\Delta) \cap M$ such that $x,y,z,w$ do not lie on a plane and with $| \langle x,y,z \rangle \cap M| = 3$. Up to adding some multiple of $w$ to $x,y,z,w$ we may assume that $x,y,z \in H_{n_{\Gamma},0}, \quad w \in H_{n_{\Gamma},-l_{\Gamma}}$ and
	
	\begin{align*}
		| \langle x,y,z,w \rangle \cap M | = 4.
	\end{align*}
	Here we implicitly used the formula of theorem (\ref{Theorem_formula_ker_ITT}). Besides apparently for $Q:= \langle x,y,z,w,y+w,z+w \rangle$ we have $|Q \cap M| = 6$ since we translate the triangle $\langle x,y,z \rangle$ by $w$ (compare corollary \ref{lemma_def_l_i_Supp_h_Gamma_i}). We apply lemma \ref{lemma_White_simplex_parallel__not_empty}. In the notation of this lemma if {\bf $w = (1,0,0)$} then $n_{\Gamma} = (-k,1,0)$ and $p=-k \geq 2$. If $w = (0,0,1)$ then $n_{\Gamma} = (0,-q,p)$ and $p = -k \geq 2$. Finally if $w= (1,p,q)$ then $n_{\Gamma} = (0,-1,0)$ and $p = -k \geq 2$.

\end{proof}

\begin{lemma} \label{lemma_White_simplex_parallel__not_empty}
	Given an $3$-dimensional empty simplex $Q_1 := \langle x=0,y,z,w \rangle$ by White's Theorem (\cite[p.34, the Terminal Lemma]{Oda88}) we may assume 
	\begin{align} \label{empty_simplex_vertices_tilde}
		\{ x \} = \begin{pmatrix} 0\\0\\0 \end{pmatrix}, \{y,z,w \} = \{ \begin{pmatrix} 1\\0\\0 \end{pmatrix}, \, \begin{pmatrix} 0\\0\\1 \end{pmatrix}, \, \begin{pmatrix} 1\\p\\q \end{pmatrix} \}
	\end{align}
	where $0 < p < q$ and $\gcd(p,q) = 1$. Assume $p \geq 2$. Then for $Q:= \langle x,y,z,w,y+w,z+w \rangle$ we get $|Q \cap M| \geq 7$. 
\end{lemma}

\begin{proof}
	If $w=(1,0,0)$:
	\begin{align*}
		Q = \langle \begin{pmatrix} 0\\0\\0 \end{pmatrix}, \begin{pmatrix} 1\\0\\0 \end{pmatrix}, \begin{pmatrix} 0\\0\\1 \end{pmatrix}, \begin{pmatrix} 1\\p\\q \end{pmatrix}, \begin{pmatrix} 1\\0\\1 \end{pmatrix}, \begin{pmatrix} 2\\p\\q \end{pmatrix} \rangle.
	\end{align*}
	Distinguish if $q$ is even and $p$ is odd let $q = r \cdot p + s$ with $r,s \in \mathbb{N}_{\geq 1}$ and $0<s<p$. Then
	\begin{align*}
		\frac{1}{p} \cdot \Big( \begin{pmatrix}
			1\\p\\q
		\end{pmatrix} + (p-s) \cdot \begin{pmatrix}
			1\\0\\1
		\end{pmatrix} + (s-1) \cdot \begin{pmatrix}
			1\\0\\0
		\end{pmatrix}  \Big) = \begin{pmatrix} 1\\1\\r+1	\end{pmatrix} \in Q \cap M
	\end{align*} 
	If $q$ is odd and $p$ is even then
	\begin{align*}
		\frac{1}{2} \cdot \Big( \begin{pmatrix}
			2\\p\\q
		\end{pmatrix} + (q-p) \cdot \begin{pmatrix}
			0\\0\\1
		\end{pmatrix} = \begin{pmatrix}
			1\\ \frac{p}{2} \\ \frac{q+1}{2}
		\end{pmatrix} \in Q \cap M
	\end{align*} 
	
	If $w=(0,0,1)$:
	\begin{align*}
		Q = \langle \begin{pmatrix} 0\\0\\0 \end{pmatrix}, \begin{pmatrix} 1\\0\\0 \end{pmatrix}, \begin{pmatrix} 0\\0\\1 \end{pmatrix}, \begin{pmatrix} 1\\p\\q \end{pmatrix}, \begin{pmatrix} 1\\0\\1 \end{pmatrix}, \begin{pmatrix} 1\\p\\q+1 \end{pmatrix} \rangle.
	\end{align*}
	If $q$ is even and $p$ is odd again $(1,1,r+1) \in Q \cap M$. If $q$ is odd and $p$ is even then $(1,\frac{p}{2}, \frac{q+1}{2}) \in Q \cap M$ as above. \leavevmode \\ \\
	Finally if $w=(1,p,q)$:
	\begin{align*}
		Q = \langle \begin{pmatrix} 0\\0\\0 \end{pmatrix}, \begin{pmatrix} 1\\0\\0 \end{pmatrix}, \begin{pmatrix} 0\\0\\1 \end{pmatrix}, \begin{pmatrix} 1\\p\\q \end{pmatrix}, \begin{pmatrix} 2\\p\\q \end{pmatrix}, \begin{pmatrix} 1\\p\\q+1 \end{pmatrix} \rangle.
	\end{align*}
	If $q$ is even and $p$ is odd let $r:= q \, \mod \, p$. There are $\lambda,\mu \geq 0$ with
	\[ \lambda + \mu = r, \quad \lambda \cdot q + \mu \cdot (q+1) = l \cdot p, \quad l \in \mathbb{Z}.  \]
	Then
	\begin{align*}
		\frac{1}{p} \cdot \big( (p-r) \cdot \begin{pmatrix} 1\\0\\0	\end{pmatrix} + \lambda \cdot \begin{pmatrix} 1\\p\\q \end{pmatrix} + \mu \cdot \begin{pmatrix}  1\\p\\q+1 \end{pmatrix} \big) = \begin{pmatrix} 1\\r\\l	\end{pmatrix} \in Q \cap M.
	\end{align*}
	If $q$ is odd and $p$ is even then $(1,\frac{p}{2}, \frac{q+1}{2}) \in Q \cap M$.
\end{proof}

	\begin{remark}
		\normalfont
	Note that we consider some known counterexamples to the ITT in (\cite{Gie23}), namely special Kanev and Todorov surfaces in toric $3$-folds, where $\Int(\Delta) \cap M = \{(0,0,0)\}$, but $\dim \, F(\Delta) = 3$ . Whereas the dimension of the kernel of $\kappa_f$ turn out to be independent of $f \in B$ the dimension of the kernel of $d \mathcal{P}_{B,f}$ in general depends on $f \in B$. This is because of the condition $g_{\Gamma}(f) \cdot x^w \in L(\Delta)$ which might very well depend on $f \in B$. 
	\end{remark}

\section{The Kodaira-Spencer map for smooth birational models} \label{section_Kod_Sp_map_smooth_bir_model}

We carry out Remark \cite[rem.9.4]{Gie22}: Replacing $\Delta$ by $\Convhull(C(\Delta) \cap M)$ we obtain $\kappa_{\mathbb{P},f} = \kappa_f$. Let $\tilde{\Sigma}$ be a smooth fan refining $\Sigma$. Let $\tilde{\mathbb{P}}$ be the toric variety to $\tilde{\Sigma}$ and $\tilde{Y} = \tilde{Y}_f$ the closure of $Z_f$ in $\tilde{\mathbb{P}}$. We 
consider the Kodaira-Spencer map
\[ \kappa_{\tilde{\mathbb{P}},f}: H^0(\tilde{Y}, N_{\tilde{Y} / \tilde{\mathbb{P}}}) \rightarrow H^1(\tilde{Y}, T_{\tilde{Y}}).  \] 
By \cite[Prop.7.1]{Bat22} $\tilde{Y} \sim_{lin} - \sum\limits_{i} \Min_{\Delta}(n_{i}) \cdot D_i$ and since $\tilde{\Sigma}$ refines $\Sigma$ and with it $\Sigma_{C(\Delta)}$ we get 
\[ H^0(\tilde{Y}, N_{\tilde{Y} / \tilde{\mathbb{P}}}) \cong L(\Delta)/ \mathbb{C} \cdot f.  \]
\begin{corollar} \label{corollary_ker_Kod_Sp_smooth_case}
	In the above situation
	\begin{align*}
		\ker(\kappa_{\tilde{\mathbb{P}},f}) \cong \Lie \, \Aut(\tilde{\mathbb{P}}).
	\end{align*}
	and for $n \geq 4$ the sequence
	\begin{align*}
		0 \rightarrow Im(\kappa_{\tilde{\mathbb{P}},f}) \rightarrow H^1(\tilde{Y}, T_{\tilde{Y}}) \rightarrow H^1(\tilde{\mathbb{P}}, T_{\tilde{\mathbb{P}}}) \rightarrow 0
	\end{align*}
	is exact
\end{corollar}

\begin{proof}
	The proof is the same as that of \cite[Thm. 6.1]{Gie22} except that we can not use the same vanishing Theorem if $\tilde{Y} + K_{\tilde{\mathbb{P}}}$ is not nef. We have to show that
	\begin{align} \label{vanishing_necessary_smooth_bir_models_ker}
		0 = H^{k}(\tilde{\mathbb{P}}, \Omega_{\tilde{\mathbb{P}}}^1 \otimes \mathcal{O}(\tilde{Y} + K_{\tilde{\mathbb{P}}})) \quad \quad k=n-2,n-1,n
	\end{align}
	remains valid.
	
	\begin{lemma} \label{toric_resolution_of_singularities}
		Let $\phi: \tilde{\mathbb{P}} \rightarrow \mathbb{P}$ be the birational toric morphism to the refinement $\tilde{\Sigma}$ of $\Sigma$. Then
		\[ 
		R^k \phi_*(K_{\tilde{\mathbb{P}}} + \tilde{Y}) = \left\{\begin{array}{ll}
			K_{\mathbb{P}} + Y & k=0 \\
			0 & k > 0
		\end{array} \right.
		\]
		Further for $D_j \subset \tilde{\mathbb{P}}$ a torus invariant prime divisor to a ray $n_j \in \tilde{\Sigma}[1]$
		\[ 
		R^k \phi_*(K_{\tilde{\mathbb{P}}} + \tilde{Y} - D_j) = \left\{\begin{array}{ll}
			K_{\mathbb{P}} + Y - D_j & k=0, \, n_j \in \Sigma[1] \\
			K_{\mathbb{P}} + Y & k = 0, \, n_j \notin \Sigma[1] \\
			0 & k > 0
		\end{array} \right.
		\]		
	\end{lemma}
	
	\begin{proof}
		We may check this locally: Let $U_{\sigma}$ be the affine toric variety to a maximal dimensional cone $\sigma$ of $\Sigma$ and assume that $n_j$ refines $\sigma$, that is
		\[ n_j = \sum\limits_{i = 1}^r \lambda_i \cdot n_i \quad n_i \in \sigma[1], \quad \lambda_i \geq 0.  \]
		Let $\delta \in \{0,1\}$. To 
		\[ K_{\tilde{\mathbb{P}}} + \tilde{Y} - \delta \cdot D_j  \]
		is associated the polytope
		\[ P= \{ x \in M_{\mathbb{R}} \mid \langle x,n_i \rangle \geq \Min_{\Delta}(n_i)+1, \, \langle x,n_j \rangle \geq \Min_{\Delta}(n_j)+1 + \delta   \}  \]
		where $n_i$ runs over all rays in $\tilde{\Sigma}$ refining the cone $\sigma$. If $n_j \notin S_F(\Delta) = \Sigma[1]$ and
		\[ \langle x,n_i \rangle \geq \Min_{\Delta}(n_i) + 1 = \Min_{F(\Delta)}(n_i)  \quad \forall \, n_i \in \sigma[1] \]
		then automatically
		\[ \langle x, n_j \rangle = \sum\limits_{i=1}^r \lambda_i \cdot \langle x, n_i \rangle \geq \Min_{F(\Delta)}(n_j) \geq \Min_{\Delta}(n_j) + 2.  \]
		Thus 
		\[ \phi_*(K_{\tilde{\mathbb{P}}} + \tilde{Y} - \delta \cdot D_j) = K_{\mathbb{P}} + Y.  \]
		The case $n_j \in \Sigma[1]$ is similar. Let $k>0$ and assume $n_j \notin \Sigma[1]$ (the other case is similar). We follow (\cite[§7,§8]{Dan78}): With respect to the action of the torus $T$
		\[ H^k(\phi^{-1}(U_{\sigma}), K_{\tilde{\mathbb{P}}} + \tilde{Y} - \delta \cdot D_j)  \]
		decomposes into the direct sum of pieces 
		\[ H^k(\phi^{-1}(U_{\sigma}), K_{\tilde{\mathbb{P}}} + \tilde{Y} - \delta \cdot D_j)_m \cong H_{Z_m}^k(\sigma, \mathbb{C})  \]
		for each $m \in M$, where
		\[ Z_m := \{ x \in \sigma \mid \langle x, m \rangle \geq \ord_{F(\Delta)}(x)  \}.  \]
		If $m \in F(\Delta)$ then $Z_m = \sigma$ and if $m \notin F(\Delta)$ then $\sigma \setminus Z_m$ is nonempty and convex. It follows from the local cohomology sequence
		\[ ... \rightarrow H^{k-1}(\sigma \setminus Z_m, \mathbb{C}) \rightarrow H_{Z_m}^k(\sigma, \mathbb{C}) \rightarrow H^k(\sigma, \mathbb{C}) = 0 \]
		that $H_{Z_m}^k(\sigma, \mathbb{C}) = 0$ for all $k \geq 1$.
	\end{proof}
	
	Use an \textit{Euler exact sequence} on $\tilde{\mathbb{P}}$ (\cite[Thm.8.1.6]{CLS11})
	\begin{align*}
		0 \rightarrow \Omega_{\tilde{\mathbb{P}}}^1 \rightarrow \bigoplus\limits_{n_i \in \tilde{\Sigma}[1]} \mathcal{O}(-D_i) \rightarrow \Pic(\tilde{\mathbb{P}}) \otimes_{\mathbb{Z}} \mathcal{O}_{\tilde{\mathbb{P}}} \rightarrow 0
	\end{align*}
	By the above Lemma and since $Y+K_{\mathbb{P}}$ is nef
	\begin{align*}
		& H^k(\tilde{\mathbb{P}}, \mathcal{O}(\tilde{Y} + K_{\tilde{\mathbb{P}}})) \cong H^k(\mathbb{P}, \mathcal{O}(Y+K_{\mathbb{P}})) = 0 \quad k > 0 \\
		& H^k(\tilde{\mathbb{P}}, \mathcal{O}(\tilde{Y} + K_{\tilde{\mathbb{P}}} - D_j)) = 0 \quad \quad k > 0
	\end{align*}
	where the last line follows by distinguishing $n_j \notin \Sigma[1]$ and $n_j \in \Sigma[1]$ and in the second case using an Euler exact sequence on $\mathbb{P}$. The vanishing (\ref{vanishing_necessary_smooth_bir_models_ker}) follows.
	
\end{proof}

\section{A realization of $\Phi_f$ on $\coker(\kappa_{\mathbb{P},f})$ by cup product and contraction} \label{section_ITT_on_coker_computation}

For the sake of notation we write $\mathbb{P}$ and $Y$ for $\tilde{\mathbb{P}}$ and $\tilde{Y}$, that is we assume $\mathbb{P}$ and $Y$ being smooth (we just work with the smooth model). Note that $D:= Y-K_{\mathbb{P}}$ is an SNC-divisor. Remember
\begin{align*}
	H^0(Y,N_{Y/\mathbb{P}}) \overset{\kappa_{\mathbb{P},f}}{\rightarrow} H^1(Y,T_Y) \rightarrow H^1(Y, T_{\mathbb{P} \vert{Y}}) \rightarrow H^1(Y,N_{Y/\mathbb{P}}) = 0
\end{align*}
The last vanishing follows since $\phi_* \tilde{Y}_f = Y_f$ similar to lemma \ref{toric_resolution_of_singularities} and $Y_f \subset \mathbb{P}$ is nef.
The following diagram shows how the cup product and contraction could be realized on $H^1(\mathbb{P}, T_{\mathbb{P}})$. 

\begin{equation*}
	\begin{tikzcd}
		H^1(\mathbb{P}, T_{\mathbb{P}}) \arrow{r} & H^1(Y,T_Y) \arrow{r}{\textrm{cup + contr.}} & \bigoplus_{p} \Hom(R_{Int,f}^{n-p+1}, R_{Int,f}^{n-p+2})  \\
		H^1(\mathbb{P}, T_{\mathbb{P}}(-\log \, Y)) \arrow{u}{\textrm{iso. up to } \kappa_f}
		&& \bigoplus_{p} \Hom(R_f^{n-p+1}, R_f^{n-p+2}) \arrow[swap]{u}{\textrm{restr.}} \arrow{d}{\cong} \\
		H^1(\mathbb{P}, T_{\mathbb{P}}(-\log \, D)) \arrow{u}{\cong} \arrow{rr}{\textrm{cup + contr.}} && \bigoplus_{p}  \Hom(H^{n-p}(\mathbb{P}, \Omega_{\mathbb{P}}^p(\log \, D)), \\
		&& H^{n-p+1}(\mathbb{P}, \Omega_{\mathbb{P}}^{p-1}(\log \, D))) 
	\end{tikzcd}
\end{equation*}

\begin{proof}
	The right vertical relations follow by \ref{definition_R_Int_f_k} and \cite[Cor.6.0.8]{Bat93}.  The first vertical homomorphism on the left is surjective since $H^1(Y,N_{Y/\mathbb{P}}) = 0$. In fact it is an isomorphism modulo $H^0(Y,N_{Y/\mathbb{P}})$. By duality
	\begin{align*}
		& H^1(\mathbb{P}, T_{\mathbb{P}}(- \log \, D)) \cong H^{n-1}(\mathbb{P}, \Omega_{\mathbb{P}}^1(\log \, D) \otimes K_{\mathbb{P}})^* 
	\end{align*}
	and similarly for $Y$ instead of $D$. Given $0 \leq D_1 \leq D = Y - K_{\mathbb{P}}$ we have the following exact sequence (\cite[2.3 Properties]{EsVi92})
	\begin{align*}
		& 0 \rightarrow \Omega_{\mathbb{P}}^1 \rightarrow \Omega_{\mathbb{P}}^1(\log \, D) \rightarrow \bigoplus\limits_{j=1}^r \mathcal{O}_{D_j} \rightarrow 0 \\
		& 0 \rightarrow \Omega_{\mathbb{P}}^p(\log \, D-D_1) \rightarrow \Omega_{\mathbb{P}}^p(\log \, D) \rightarrow \Omega_{D_1}^{p-1}(\log \, D-D_1) \rightarrow 0
	\end{align*}

	Let $-K_{\mathbb{P}} = \sum\limits_{i} D_i$ Then
	\begin{align*}
		&... \rightarrow H^{n-1}(\mathbb{P}, \Omega_{\mathbb{P}}^1(\log \, Y) \otimes K_{\mathbb{P}}) \rightarrow H^{n-1}(\mathbb{P}, \Omega_{\mathbb{P}}^1(\log \, D) \otimes K_{\mathbb{P}}) \\
		&\rightarrow \bigoplus_{i} H^{n-1}(D_i, \mathcal{O}_{D_i}(Y + K_{\mathbb{P}} + C_i)_{\vert{D_i}}) \rightarrow ...
	\end{align*}
	with $0 \leq C_i < -K_{\mathbb{P}}$. Further we have an exact sequence
	\begin{align*}
		0 \rightarrow \mathcal{O}_{D_i}(Y+K_{\mathbb{P}} + C_i)_{\vert{D_i}} \rightarrow \mathcal{O}_{D_i}(Y) \rightarrow \bigoplus_{j} \mathcal{O}_{D_i \cap D_j}(Y_{\vert{D_i \cap D_j}})
	\end{align*}
	Then $\mathcal{O}(Y_{\vert{D_i}})$ and $\mathcal{O}(Y_{\vert{D_i \cap D_j}})$ are nef divisors on toric varieties. Thus for $n \geq 4$
	\begin{align*}
		H^{n-k}(D_i, \mathcal{O}_{D_i}(Y + K_{\mathbb{P}} + C_i)_{\vert{D_i}}) = 0 \quad k=1,2.
	\end{align*}

\end{proof}

\section{A duality homomorphism and the ITT}
\label{section_ITT_on_coker_computation_2}

Let $D_1 := Y$. Then the first sequence from the proof above gives
\begin{align*}
	0 \rightarrow \Omega_{\mathbb{P}}^1(\log \, -K_{\mathbb{P}}) \rightarrow \Omega_{\mathbb{P}}^1(\log \, D) \rightarrow \mathcal{O}_Y(-K_{\mathbb{P} \vert{Y}}) \rightarrow 0.
\end{align*}
Tensor this sequence with $K_{\mathbb{P}}$:
\begin{align*}
	H^{n-1}(\mathbb{P}, \Omega_{\mathbb{P}}^1(\log \, D) \otimes K_{\mathbb{P}}) \subset \underbrace{H^{n-1}(\mathbb{P}, \Omega_{\mathbb{P}}^1(\log \, -K_{\mathbb{P}}) \otimes K_{\mathbb{P}})}_{\cong H^1(\mathbb{P}, T_{\mathbb{P}}(-\log \, -K_{\mathbb{P}}))^*=0} \oplus H^{n-1}(Y, \mathcal{O}_Y)
\end{align*}
since $T_{\mathbb{P}}(-\log \, -K_{\mathbb{P}}) \cong \mathcal{O}_{\mathbb{P}} \otimes_{\mathbb{Z}} N$. $H^{n-1}(Y,\mathcal{O}_Y)^* \cong H^0(Y,\mathcal{O}_Y(K_Y))$. Besides
\begin{align*}
	H^{0}(Y,\mathcal{O}_Y(K_Y)) \overset{\Phi_f^{n-1}}{\rightarrow} \Hom \Big( \mathbb{C}^{\vee}, H^{n-1}(Y,\mathcal{O}_Y) \Big) \cong \Hom \Big( H^0(Y,\mathcal{O}_Y(K_Y), \mathbb{C} \Big)
\end{align*}
given by cup product and contraction defines a perfect pairing, where $\mathbb{C} \cong \mathbb{C} \cdot x^{(0,...,0)} \subset R_{Int,f}^2 \approx H^1(Y, \Omega_Y^{n-2})$.

\begin{theorem}
	Let $\Delta \subset M_{\mathbb{R}}$ be an $n$-dimensional lattice polytope with $n \geq 4$ and $\Int(\Delta) \cap M \nsubseteq \textrm{plane}$ and $f \in L(\Delta)$ is nondegenerate. Then given a functorial toric resolution of singularities $\tilde{\mathbb{P}} \rightarrow \mathbb{P}$ (such a resolution exists by \cite[Rem.2.3.10]{Wlo26}) the smooth birational model $Y = \tilde{Y}_f$ fulfills the infinitesimal Torelli theorem (ITT), that is
	\begin{align*}
		\Phi_f = \bigoplus_p \Phi_f^p: H^1(Y,T_Y) \rightarrow \bigoplus_{p} \Hom \Big( H^{n-1-p}(Y, \Omega_Y^p), H^{n-p}(Y, \Omega_Y^{p-1}) \Big)
	\end{align*}
	is injective. 
	
	\begin{proof}
		Replacing $\Delta$ by $\Convhull(C(\Delta) \cap M)$ we obtain $\kappa_{\mathbb{P},f} = \kappa_f$ (\cite[Prop.9.2]{Gie22}). By Corollary \ref{corollary_ker_Kod_Sp_smooth_case} given $n \geq 4$ we have 
		\begin{align*}
			\ker(\kappa_{\tilde{\mathbb{P}},f}) \cong H^0(\tilde{\mathbb{P}}, T_{\tilde{\mathbb{P}}}) \cong H^0(\mathbb{P}, T_{\mathbb{P}})
		\end{align*}
		(see \cite[Cor.4.7]{GKK10}) since $\tilde{\mathbb{P}} \rightarrow \mathbb{P}$ is a functorial resolution of singularities. Further we obtain an exact sequence
		\begin{align*}
			0 \rightarrow Im(\kappa_{\mathbb{P},f}) \rightarrow H^1(Y, T_{Y}) \rightarrow H^1(\mathbb{P}, T_{\mathbb{P}}) \rightarrow 0.
		\end{align*}
		We may identify 
		\begin{align*}
			Im \, \Phi_f^{n-1}: \Hom(R_{Int,f}^{n-1}, R_{Int,f}^n) \cong \Hom(R_{Int,f}^1, R_{Int,f}^2)
		\end{align*}
		where $R_{Int,f}^1 \cong L^*(\Delta)$. Next given the assumptions $\Int(\Delta) \cap M \nsubseteq \{ \textrm{plane}\}$ we get
		\begin{align*}
			\ker(\Phi_{f \vert{Im \, \kappa_f}}^{n-1}) \overset{\textrm{duality}}{=} \ker(\Phi_{f \vert{Im \, \kappa_f}}^{1}) = \{0\}
		\end{align*}
		(Corollary \ref{corollary_inf_Tor_thm}). Besides $\ker(\Phi_{f \vert{\coker \, \kappa_f}}) \subset \ker( \Phi_{f \vert{\coker \, \kappa_f}}^{n-1}) = \{0\}$ by section \ref{section_ITT_on_coker_computation} and \ref{section_ITT_on_coker_computation_2}. By assumption $l^*(\Delta) \geq 2$ and $\Phi_{f \vert{\coker \, \kappa_f}}^{n-1}$ defines a perfect pairing on $L^*(\Delta) \cong H^0(Y,\mathcal{O}_Y(K_Y))$ with $\Phi_{f \vert{\coker \, \kappa_f}}^{n-1}(L^*(\Delta)) = \mathbb{C} \cdot x^{(0,...,0)}$. Thus certainly
		\begin{align*}
			Im \, \Phi_{f \vert{Im \, \kappa_{f}}}^{n-1} \cap \, Im \, \Phi_{f \vert{\coker \, \kappa_f}}^{n-1} = \{0\}
		\end{align*}
		and it follows $\ker(\Phi_f) = \{0\}$ proving the assertion.
	\end{proof}
	
	\begin{remark}
		\normalfont
		Assuming that $\coker(\kappa_{\mathbb{P},f})$ does not induce difficulties for $n=3$ we would get: Given an $n$-dimensional lattice polytope $\Delta$, $n \geq 2$, with $\Int(\Delta) \cap M \nsubseteq \textrm{ plane}$ and $f$ generic then $Y$ fulfills the ITT. This excludes for $n=2$ the known case of hyperelliptic curves. For $n=3$ we found explicit counterexamples where $\Int(\Delta)\cap M = \{(0,0,0)\}$ (Kanev and Todorov surfaces), and $\Int(\Delta) \cap M \subseteq \textrm{ line}$ (elliptic surfaces).
	\end{remark}
	
\end{theorem}

\end{document}